\documentclass[final]{siamltex}

\usepackage{times}

\usepackage{amsmath,amssymb,dsfont,graphicx,eucal,color,comment,epsfig,overpic}

\newcommand{\T}{\textrm}

\newcommand{\be}{\begin{equation}}
\newcommand{\ee}{\end{equation}}
\newcommand{\bea}{\begin{eqnarray}}
\newcommand{\eea}{\end{eqnarray}}
\newcommand{\beaa}{\begin{eqnarray*}}
\newcommand{\eeaa}{\end{eqnarray*}}

\newtheorem{remark}[theorem]{Remark}

\newcommand{\mbf}[1]{\mbox{\boldmath$\rm{#1}$}}
\newcommand{\smbf}[1]{\mbox{\boldmath \scriptsize$\rm{#1}$}}

\newcommand{\ddd}{\text{\rm D}}
\newcommand{\dn}{\text{\rm N}}
\newcommand{\el}{ \kappa \in \mathcal{T}}

\newcommand{\dint}{\text{\rm int}}

\newcommand{\ha}{\frac{1}{2}}

\newcommand{\su}{\sum_{\kappa\in\mathcal{T}}}\newcommand{\ud}{\mathrm{d}}
\newcommand{\ndg}[1]{| \kern -.25mm \|{#1}| \kern -.25mm \|}
\newcommand{\ltwo}[2]{\|{#1}\|_{#2}}
\newcommand{\ltwoin}[2]{\langle{#1},{#2}\rangle}
\newcommand{\qed}{\vspace{-.2cm} \\* \mbox{} \hfill $\Box$ }

\newcommand{\ex}{{\rm e}}
\newcommand{\texte}[1]{\quad \text{{#1}}\quad}
\newcommand{\honen}{\mathbb{H}^1}
\newcommand{\Hn}{\mathbb{H}}
\newcommand{\FEspace}{{\rm V}_h}
\newcommand{\FESpace}{\mathbb{V}_h}
\newcommand{\SSpace}{\mathbb{S}}

\newcommand{\Gint}{\Gamma_{\dint}}

\newcommand{\inter}{{\mathcal {I}}}
\newcommand{\Gatr}{\Gamma_{\inter}}

\newcommand{\bu}{{\mbf{u}}}
\newcommand{\bv}{{\mbf{v}}}
\newcommand{\bw}{{\mbf{w}}}

\newcommand{\bz}{{\mbf{z}}}

\newcommand{\bxi}{\mbf{\xi}}

\newcommand{\boeta}{\mbf{\eta}}
\newcommand{\brho}{\mbf{\rho}}
\newcommand{\sbrho}{\smbf{\rho}}
\newcommand{\btheta}{\mbf{\theta}}
\newcommand{\pjump}{\tilde{\mbf{p}}}
\newcommand{\bx}{{\mbf{x}}}

\newcommand{\tth}{\mbf{\tt h}}

\newcommand{\mean}[1]{ \{#1\} }
\newcommand{\jump}[1]{  [#1]  }
\newcommand{\jumptwo}[1]{[\![#1]\!]}
\newcommand{\weight}{\upsilon}
\newcommand{\Weight}{\Upsilon}
\newcommand{\friction}{r}
\newcommand{\Friction}{{\rm R}}

\DeclareMathOperator{\diam}{diam}

\DeclareMathOperator*{\esssup}{ess\,sup}

\usepackage{morefloats}
\usepackage{xifthen}
\usepackage{color}
\usepackage{pifont}
\definecolor{darkgreen}{RGB}{0,150,0}

\newboolean{draft}
\setboolean{draft}{true}

\ifthenelse{\boolean{draft}}{}{}
\ifthenelse{\boolean{draft}}{}{}
\ifthenelse{\boolean{draft}}{}{}
\ifthenelse{\boolean{draft}}{\newcommand{\caa}[1]{\marginpar{\scriptsize{\color{red} #1}}}}{\newcommand{\caa}[1]{}}
\ifthenelse{\boolean{draft}}{}{}
\ifthenelse{\boolean{draft}}{}{}
\ifthenelse{\boolean{draft}}{}{}
\ifthenelse{\boolean{draft}}{}{}
\ifthenelse{\boolean{draft}}{}{}

\title{Discontinuous Galerkin Methods\\ for Mass Transfer \\ through Semi-Permeable Membranes}

\author{Andrea Cangiani
\thanks{Department of Mathematics, University of Leicester, University Road,
Leicester LE1 7RH, United Kingdom,
e-mail: {\tt andrea.cangiani@le.ac.uk}}
\and
Emmanuil H. Georgoulis
\thanks{Department of Mathematics, University of Leicester, University Road,
Leicester LE1 7RH, United Kingdom,
e-mail: {\tt Emmanuil.Georgoulis@le.ac.uk}}
\and
Max Jensen
\thanks{Department of Mathematical Sciences, University of Durham, Science Laboratories, South Road, Durham DH1 3LE, United Kingdom,
e-mail: {\tt m.p.j.jensen@durham.ac.uk}}
}

\begin{document}

\maketitle
\begin{abstract}
A discontinuous Galerkin (dG) method for the numerical solution of initial/boundary value multi-compartment partial differential equation (PDE) models, interconnected with interface conditions, is presented and analysed.  The study of interface problems is motivated by models of mass transfer of solutes through semi-permeable membranes. More specifically, a model problem consisting of a system of semilinear parabolic advection-diffusion-reaction partial differential equations in each compartment, equipped with respective initial and boundary conditions, is considered. Nonlinear interface conditions modelling selective permeability, congestion and partial reflection are applied to the compartment interfaces. An interior penalty dG method is presented for this problem and it is analysed in the space-discrete setting. The a priori analysis shows that the method yields optimal a priori bounds, provided the exact solution is sufficiently smooth. Numerical experiments indicate agreement with the theoretical bounds and highlight the stability of the numerical method in the advection-dominated regime.

\end{abstract}
\section{Introduction}

Models of mass transfer of substances (solutes) through semi-perme\-able membranes appear in various contexts, such as biomedical and chemical engineering applications~\cite{Friedman}. Examples include the modelling of electrokinetic flows, solute dynamics across arterial walls, and cellular signal transduction (see, e.g., \cite{Brera2010,Zunino,Cangiani_Natalini_2010} and the references therein).

This work is concerned with the development and analysis of numerical methods for a class of continuum models for mass transfer based on initial/boundary value multi-compart\-ment partial differential equation (PDE) problems, closed by nonlinear interface conditions. The interface conditions considered are the Kedem-Katchalsky (KK) equations, which represent an established model for the mass transfer mechanisms~\cite{KK58,Katchalsky}. More specifically, we consider a generic model problem consisting of a system of semilinear advection-diffusion-reaction parabolic PDE problems in multi-compartment configurations, coupled with nonlinear interface KK-type conditions. The focus is to address some challenges in the numerical solution of these models, such as the treatment of nonlinearities due to both the interface modelling and the nonlinear reactions, the discontinuity of the state variables across the interface, as well as the development of stable numerical methods in the advection-dominated regime.

Numerical methods for mass transfer problems based on conforming finite elements have been developed for the solution of solute dynamics across arterial walls; see~\cite{Zunino,Quarte02,Quarte05} and the references therein for more details. Some existence results for the purely diffusing interface problem without forcing, coupled with KK-type interface conditions, along with some numerical experiments are given in~\cite{Calabro06}. Further, numerical approaches to the treatment of interface conditions for PDE problems, resulting to globally continuous solutions can be found, e.g., in~\cite{babuska,AkrivisDougalis91,MR1622502,MR2013126,Melenk_interface10}.

Discontinuous Galerkin (dG) methods (see, e.g., \cite{nodal_dg_book,riviere_book,dipietro_ern_book} and the references therein) 
and mortar methods (see, e.g., \cite{mortar_survey} for a survey) have also been proposed for the treatment of coupled systems via interface conditions in various contexts~\cite{Stenberg98,girault_riviere,electroporation,ern_moz,MR2208566,kan_riv,MR1975269,Hansbo05}. Also, the advantages of dG methods for interfacing different numerical methods (numerical interfaces) have been identified~\cite{per_scho,MR1956028}, as well as their use on transmission-type/high-contrast problems, yielding continuous solutions across the transmission interface, has been investigated~\cite{MR2002258,MR2257119,MR2491426,MR2837483}.

Here, we consider a dG method of interior penalty type for the solution of the semilinear parabolic advection-diffusion-reaction PDE system coupled with nonlinear interface conditions of KK-type across the subdomains. The use of dG is motivated partly by the observation that the interface conditions, yielding \emph{discontinuous} solutions across the interface, can be imposed by modifying the interior penalty dG numerical fluxes. Another important factor for employing a dG method is the desired stability property of the numerical method in the advection-dominated regime.  A priori bounds for the  proposed spatially discrete dG method in both the $L^\infty(L^2)$- and $L^2(H^1)$-type norms are presented for a range of reaction fields, under the simplifying assumption that the finite element mesh is aligned with the subdomain interfaces. 

A priori error bounds for interior penalty dG methods for parabolic problems have been considered in various settings (see, e.g.,~\cite{riviere_book} for an exposition and the more recent~\cite{Dolej08}).  DG methods for semilinear parabolic spatially self-adjoint problems with locally Lipschitz continuous nonlinearity have been analysed in~\cite{Lasis07}.
In the present analysis, advection terms are included and systems of equations are considered. In the presence of advection, the analysis of the symmetric interior penalty dG method in~\cite{Lasis07} would require the assumption of quasi-uniformity of the mesh. To avoid this assumption, a different continuation argument is employed in the derivation of the a priori bounds presented here, at the expense of a stricter growth condition on the nonlinearity of the forcing term. This continuation argument is inspired by the derivation of a posteriori bounds for semilinear parabolic phase-field models~\cite{nochetto,bartels}.  The nonlinear interface terms are tackled using a non-standard elliptic projection which is inspired by a classical construction of Douglas and Dupont~\cite{douglas_dupont} for the treatment of nonlinear boundary conditions. 

The paper is organised as follows. In Section~\ref{prelim}, the PDE model along with a short derivation of the nonlinear interface conditions is presented. Section~\ref{dg_heat} is devoted to the description of the dG method proposed for the advection-diffusion part of the spatial operator incorporating the nonlinear interface conditions. Section~\ref{err_ell_problem} contains error estimates for the new elliptic projection, which is, in turn, utilised in the subsequent a priori error analysis presented in Section~\ref{err_par_prob}. Section~\ref{sec:numex} contains numerical experiments highlighting the stability and the optimal rate of convergence of the proposed method in practice. Finally, some conclusions are offered in Section~\ref{conclusions}.

\section{Interface modelling and governing PDEs}\label{prelim}

We shall consider systems of parabolic semilinear PDEs describing the flux of solutes around and through a semi-permeable membrane. The membrane is modelled as an internal boundary equipped with nonlinear interface conditions which are described in the following section.

\subsection{Interface modelling}\label{sec:tras}

We outline the Kedem-Katchalsky (KK) equations modelling solutes flow across semi-permeable membranes. The KK equations have been introduced in~\cite{KK58}; we refer to~\cite{Teorell53} for earlier works and to~\cite{Friedman} for a thorough exposition.

It is assumed that the membrane separates two compartments $\Omega^1$ and $\Omega^2$ filled with a free fluid which is called the solvent and that the membrane characteristics are uniform in space and time. The KK equations specify the dependence of the solutes and solvent fluxes across the membrane in terms of two driving forces, namely the hydrostatic and osmotic pressure jumps. In the case of a single solute, the solvent and solute fluxes from $\Omega^1$ to $\Omega^2$ normal to the membrane walls are given by, respectively,
\begin{eqnarray}\label{eq:KK1}
    \vspace{2 mm}
    J_v&=& L_P (\delta {\rm p}-\sigma_d \, \delta\pi),\\
    \label{eq:KK2}
    J_s&=&\omega \, \delta\pi+J_v(1-\sigma_f)\bar{u},
\end{eqnarray}
in terms of the hydrostatic and osmotic pressure jumps $\delta {\rm p}={\rm p}^1-{\rm p}^2$ and $\delta \pi=RT\delta u$, with $R$ being the ideal gas constant, $T$ denoting temperature, and $\delta u=u^1-u^2$ the solute concentration jump across the membrane.  Here, $\bar{u}$ represents the average concentration of the solute across the membrane.  The above constitutive laws are characterised by the phenomenological coefficients of filtration $L_P$, reflection $\sigma_d$ and $\sigma_f$, and permeation $\omega$. These coefficients may depend on the concentration while they are assumed to be constant with respect to both the position along the membrane and time.  

Equation~\eqref{eq:KK1}, which is known as Starling's law of filtration, shows that the solvent flow is affected by the osmotic flow of the solute. This observation introduces a nonlinearity in the transport term of the solute flux in equation~\eqref{eq:KK2}. Indeed, substituting $J_v$ into~\eqref{eq:KK2} we get the final model for the solute flux:
\[
J_s=\omega \, R \, T \, \delta u + L_P ( \delta {\rm p} - \sigma_d R T \delta u)(1-\sigma_f) \bar{u}
=p(u^1,u^2) \, \delta u- \friction    \bar{u}(\mbf{b}\cdot \mbf{n})|_{\Omega^2}.
\]
In this last expression, we have collected the diffusive part of the flux by introducing the nonlinear {\em permeability} function $p(u^1,u^2)$, and written the advective transport in terms of the friction coefficient $\friction\in [0,1]$ and the normal component of the transport field yielded by the hydrostatic pressure, denoted by $\mbf{b}\cdot \mbf{n}$. 

In presence of $n$ solutes, it is usually assumed that there is no direct coupling between the solutes' fluxes~\cite{Friedman}. Under this simplifying assumption the KK equations describing the flux of the solvent and $n$ solutes read
\begin{eqnarray}\label{eq:KK1_many}
    J_v&=&L_P \big(\delta {\rm p}-\sum_{j=1}^n\sigma_{d,j} \, \delta\pi_j\big),\\
    \label{eq:KK2_many}
    J_{s,i} &=&\omega \, \delta\pi_i+J_v(1-\sigma_f)\bar{u}_i,\qquad i=1,\dots ,n,
\end{eqnarray}
with $J_{s,i}$, $\delta \pi_i$, and $\bar{u}_i$ representing, respectively, the flux, osmotic pressure, and average concentration of the $i$-th solute.
Proceeding as in the case of one solute we get the following expression for the solute fluxes:
\begin{eqnarray}
\nonumber
J_{s,i}&=&
p_i(u_i^1,u_i^2)
 \delta u_i-\sum_{j=1}^n \tilde{p}_{i,j}(\mbf{u}^1,\mbf{u}^2) \, \delta u_j
-\friction_i \bar{u}_i (\mbf{b}\cdot \mbf{n})|_{\Omega^2}
\\
\label{eq:KK2_final}
&=&
\mbf{p}_i(\mbf{u}^1,\mbf{u}^2)\cdot \delta\mbf{u}
-\friction_i \bar{u}_i (\mbf{b}\cdot \mbf{n})|_{\Omega^2},
\qquad i=1,\dots,n,
\end{eqnarray}
with $\delta\mbf{u}=\mbf{u}^1-\mbf{u}^2$, $\mbf{u}^j=(u_1^j,\dots,u_n^j)^{{\mathsf T}}$, $j=1,2$, where $u_i^j$ denotes the concentration of the $i$-th solute in the $j$-th compartment,  $p_i$ the corresponding permeability, and $\tilde{p}_{i,j}$ the cross-coefficients expressing the crowding effect inside the membrane. For the second equality, we have collected in the vector {\em permeability} function $\mbf{p}_i(\mbf{u}^1,\mbf{u}^2)$ all diffusion terms. 

It remains to fix a model for the average concentrations inside the membrane.  In the case of relatively thick membranes~\cite{Zunino,Quarte05}, it is appropriate to consider a mechanical approach describing solutes flow within the membrane through an advection-diffusion model. This approach yields a weighted arithmetic average of the concentration at the membrane's faces. Thus, denoting by $u^1_i$ and $u^2_i$, $i=1,\dots, n$, such concentrations, we get 
\begin{equation} \label{eq:tiltedaverage}
\bar{u}_i=\delta_w u_i:=v^1_iu^1_i+v^2_iu^2_i,
\end{equation}
for some given weights $v_i=(v^1_i,v^2_i)$ with $v^1_i+v^2_i=1$. These can be expressed in terms of the internal P\`eclet number of the membrane advection-diffusion model and are so that the upwind value dominates~\cite{Zunino}, {\em cf.} the conditions given 
below in~\eqref{eq:weights_conditions}.

In what follows, the fluxes given by~\eqref{eq:KK2_final} together with the model \eqref{eq:tiltedaverage} for the average concentration inside the membrane are used to close the PDE problem with appropriate interface conditions.


\subsection{Notation}\label{not}

We denote by $L^p(\omega)$, $1\le p\le +\infty$, the standard Lebesgue spaces, $\omega\subset\mathbb{R}^d$, $d \in \{ 2,3\}$, with corresponding norms $\|\cdot\|_{p,\omega}$; if $\omega=\Omega$, we shall write instead $\|\cdot\|_{p}$. Also the norm of $L^2(\omega)$ will be denoted by $\ltwo{\cdot}{\omega}$ and if $\omega=\Omega$ by $\ltwo{\cdot}{}$ for brevity; by $\langle \cdot,\cdot\rangle$ we write the standard $L^2$-inner product on $\Omega$; when the arguments are vectors of $L^2$-functions, the $L^2$-inner product is modified in the standard fashion.  We denote by $H^s(\omega)$ the standard Hilbertian Sobolev space of index $s\in\mathbb{R}$ of real-valued functions defined on $\omega\subset\mathbb{R}^d$; in particular $H^1_0(\omega)$ signifies the space of functions in $H^1(\omega)$ whose traces onto the boundary $\partial\omega$ vanish.  For $1\le p\le +\infty$, we denote the standard Bochner spaces $L^p(0,T;X)$, with $X$ being a Banach space with norm $\|\cdot\|_X$. Finally, we denote by $C(0,T;X)$ the space of continuous functions $v: [0,T]\to X$ with norm $\|v\|_{C(0,T;X)}:=\max_{0\le t\le T}\|v(t)\|_{X}<+\infty$.

\begin{figure}
\centering
\setlength{\unitlength}{2.5cm}
\begin{picture}(0,1.2)
\put(-1,0){\line(0,1){1}}
\put(-1,0){\line(1,0){2}}
\put(-1,1){\line(1,0){2}}
\put(1,0){\line(0,1){1}}
\put(0,0){\line(0,1){1}}
\put(-.5,.5){\makebox(0,0){$\Omega^1$}}
\put(.5,.5){\makebox(0,0){$\Omega^2$}}
\put(.1,.5){\makebox(0,0){{{$\Gatr$}}}}
\end{picture}
\caption{The domain of solution $\Omega$ is subdivided into two subdomains $\Omega^1$, $\Omega^2$. The interface boundary is defined as 
$\Gatr:=(\partial {\Omega}^1\cap\partial{\Omega}^2)\backslash \partial\Omega$.}
\label{fig:two_compartments}
\end{figure}
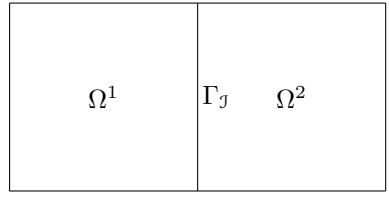

Let $\Omega$ be a bounded open domain with Lipschitz boundary in $\mathbb{R}^d$, and let $\partial\Omega$ be the boundary of $\Omega$. The domain $\Omega$ is subdivided into two subdomains $\Omega^1$ and $\Omega^2$, such that $\Omega=\Omega^1\cup\Omega^2\cup\Gatr$, where $\Gatr:=(\partial {\Omega}^1\cap\partial{\Omega}^2)\backslash \partial\Omega$, see Figure~\ref{fig:two_compartments}. For $i=1,2$, we assume that $\Omega^i$ has Lipschitz boundary and that $\partial \Omega^i\cap\partial\Omega$ has positive $(d-1)$-dimensional (Hausdorff) measure. We define $\Hn^s:=[H^s(\Omega^1\cup\Omega^2)]^{n}$, $s \in \mathbb{R}$.

We shall employ the following notational convention: vectors are indicated with lower case bold symbols, $n\times n$ diagonal matrices with upper case (non-bold) symbols, and $n\times d$ tensors with upper case bold symbols.

The gradient $\nabla \bv$ of a vector function $\bv:\Omega^1\cup\Omega^2\to\mathbb{R}^n$ in $\honen$ is a mapping $\Omega^1\cup\Omega^2\to\mathbb{R}^{n\times d}$ gained from componentwise application of the gradient operation: $\nabla \bv:=(\nabla v_1,\dots,\nabla v_n)^{{\mathsf T}}$. Similarly the divergence $\nabla \cdot\mbf{Q}$ of the tensor-valued function $\mbf{Q}:\Omega^1\cup\Omega^2\to\mathbb{R}^{n\times d}$ is $\nabla \cdot\mbf{Q}:=(\nabla \cdot Q_1,\dots,\nabla \cdot Q_n)^{{\mathsf T}}$ where the $Q_i$ are rows of $\mbf Q$.

\subsection{Model Problem}\label{mod}

For a time interval $[0,T]$, $T>0$, and for 
\[
\bu\in L^2(0,T;\honen), \qquad \bu_t\in L^2(0,T;\Hn^{-1}),
\]
with $\bu:=(u_1,\dots,u_n)^{{\mathsf T}}$, we consider the system of semilinear parabolic equations
\begin{equation}\label{pde}
 \bu_t - \nabla\cdot (A\nabla \bu - U\mbf{B})+\mbf{F}(\bu)=\mbf{0}\quad\text{in } (0,T]\times(\Omega^1\cup\Omega^2),
\end{equation}
with $U$ denoting the diagonal matrix $U=\diag(u_1,\dots,u_n)$.  Here, $\mbf{B}$ is an $n\times d$ tensor field with rows $B_i\in C^1(0,T;W^{1,\infty}({\Omega}\backslash\Gatr)^d\cap W^{\infty}({\rm div},\Omega))$, $i=1,\dots, n$, and $A\in [C([0,T] \times \Omega^1\cup\Omega^2)]^{n\times n}$ diagonal, with $A=\diag(a_1,a_2,\dots,a_n)$, where $a_i:[0,T]\times \Omega^1\cup\Omega^2\to\mathbb{R}$, $i=1,\dots, n$.  We assume that there exists a constant $\alpha_{\min}>0$ of uniform parabolicity such that $a_i(t,x)\ge\alpha_{\min}$ for all $i=1,\dots,n$ and $(t,x) \in [0,T] \times \Omega$. For simplicity, we also require that the matrix $1/2\diag(\nabla\cdot \mbf{B})$ is positive semi-definite. Finally, $\mbf{F}: \mathbb{R}^n\to\mathbb{R}^n$ is a vector field satisfying the growth condition 
\be\label{growth} 
|\mbf{F}(\bw)-\mbf{F}(\bv)|\le C(1+|\bw|+|\bv|)^{\gamma}|\bw-\bv|,
\ee
for $\bw, \bv\in\mathbb{R}^n$ and $\gamma \ge 0$ constant, where $|\cdot|$ denotes the Euclidean distance on $\mathbb{R}^n$. The admissible values of the constant $\gamma$ will be discussed in detail at various instances in the text.

We impose the initial condition
\be\label{modelbc}
  \bu(0,x)=\bu_0(x)\text{ on }\{0\}\times\Omega,
\ee
for $\bu_0\in [L^2(\Omega)]^n$. On $[0,T]\times\partial\Omega$, we impose mixed Dirichlet and Neumann boundary conditions as follows. Let $i$ be an index running over the set $\{1,\ldots,n\}$. For each $i$ the boundary $\partial\Omega$ is split into $\partial\Omega=\Gamma^{\ddd}_i\cup\Gamma^{\dn}_i$, with $\Gamma^{\ddd}_i$ being of positive $(d-1)$-dimensional (Hausdorff) measure. Further, we subdivide $\partial\Omega=\partial\Omega_i^-\cup\partial\Omega_i^+$, where $\partial\Omega_i^-:=\{\mbf{x}\in\partial\Omega: (B_i\mbf{n})(\mbf{x})<0\}$ and $\partial\Omega_i^+=\partial\Omega\backslash\partial\Omega_i^-$ are the \emph{inflow} and \emph{outflow} parts of the boundary $\partial\Omega$ for the $i$-th equation. Finally, we assign
\begin{equation}\label{bc}
\begin{aligned}
u_i &= g^{\ddd}_i \quad\text{on}\ \Gamma^{\ddd}_i, \\
a_i\nabla u_i\cdot \mbf{n} &= g^{\dn}_i \quad\text{on}\ \Gamma^{\dn}_i\cap\partial\Omega_i^+,\\
\big(a_i\nabla u_i - B_i^{\mathsf T}u_i\big)\cdot \mbf{n} &= g^{\dn}_i \quad\text{on}\ \Gamma^{\dn}_i\cap\partial\Omega_i^-,
\end{aligned}
\end{equation}
for Dirichlet and Neumann data $g^{\ddd}_i\in H^{1/2}(\Gamma^{\ddd}_i)$, $g^{\dn}_i\in L^2(\Gamma^{\dn}_i)$, respectively; here and in what follows $\mbf{n}$ denotes the unit outward normal vector to $\partial\Omega$.  We denote by $\chi_i^-:\partial\Omega_i^-\to\mathbb{R}$ the characteristic function of $\partial\Omega_i^-$.  Upon defining $\mathcal{X}^-:=\diag(\chi_1^{-},\dots,\chi_n^{-})$,  and  $\mathcal{X}^+:=I-\mathcal{X}^-$, we can write (\ref{bc}) collectively as
\be\label{bcvec}
\bu = \mbf{g}_{\ddd} \quad\text{on}\ \Gamma_{\ddd}\quad\text{and}\quad \big(A\nabla \bu-\mathcal{X}^-U\mbf{B}\big)\mbf{n} = \mbf{g}_{\dn} \quad\text{on}\ \Gamma_{\dn},
\ee
where
$\mbf{g}_{\ddd}:=(g^{\ddd}_1,\dots,g^{\ddd}_n)^{{\mathsf T}}$ and $\mbf{g}_{\dn}:=(g^{\dn}_1,\dots,g^{\dn}_n)^{{\mathsf T}}$.

The model problem is completed imposing, across $\Gatr$, the fluxes described in Section~\ref{sec:tras}. In view of~\eqref{eq:KK2_final} and \eqref{eq:tiltedaverage}, we define the friction coefficients and weights $\friction_i,\weight_i^{1,2}:\Gatr\rightarrow [0,1]$ and the {\em permeabilities} $\mbf{p}_i:\mathbb{R}^{2n}\to\mathbb{R}^n$, $i=1,\dots,n$, as functions of the traces of $\bu$ from both sides of the interface.  The interface conditions read
\[
\begin{aligned}
(a_i\nabla u_i-u_iB_i^{\mathsf T})\cdot \mbf{n}|_{\Omega^1} &= \mbf{p}_i(\bu^1,\bu^2) \cdot
(\bu^2-\bu^1)
-\friction_i(\weight_i^1 u_i^1+\weight_i^2u_i^2)(B_i\mbf{n})|_{\Omega^1},\quad\text{on}\ \Gatr,\\
(a_i\nabla u_i-u_iB_i^{\mathsf T})\cdot \mbf{n}|_{\Omega^2} &= \mbf{p}_i(\bu^1,\bu^2)\cdot
(\bu^1-\bu^2)
-\friction_i(\weight_i^1 u_i^1+\weight_i^2u_i^2)(B_i
\mbf{n})|_{\Omega^2},\quad\text{on}\ \Gatr,
\end{aligned}
\]
where $\bu^j:=\bu|_{\bar{\Omega}_j\cap\Gatr}$, $j=1,2$. Introducing,
$
\Weight^{j} =\T{diag}(\weight^{j}_1,\dots,\weight^{j}_n), \quad  j \in \{1,2\}$,
$
\Friction =\T{diag}(\friction_1,\dots,\friction_n)$, and
$\mbf{P}(\bu) =(\mbf{p}_1(\bu^1,\bu^2),\dots,\mbf{p}_n(\bu^1,\bu^2))^{{\mathsf T}}$,
the interface conditions can be written in vector notation as
\begin{equation}\label{transmission}
\begin{aligned}
(A\nabla \bu-U\mbf{B}) \, \mbf{n}|_{\Omega^1} &= \mbf{P}(\bu) (\bu^2-\bu^1)
-\Friction(\Weight^1 U^1+\Weight^2U^2)(\mbf{B}\mbf{n})|_{\Omega^1},\quad\text{on}\ \Gatr,\\
(A\nabla \bu-U\mbf{B}) \, \mbf{n}|_{\Omega^2} &= \mbf{P}(\bu) (\bu^1-\bu^2)
-\Friction(\Weight^1 U^1+\Weight^2U^2)(\mbf{B}\mbf{n})|_{\Omega^2},\quad\text{on}\ \Gatr.
\end{aligned}
\end{equation}

We make the following assumptions on the weights and permeabilities in accordance with Section~\ref{sec:tras}. For every $i=1,\dots,n$, the weights $\weight_i^{1,2}$ satisfy, for any  $\mbf{x}\in\Gatr$,
\begin{equation}\label{eq:weights_conditions}
\weight_i^1(\mbf{x})+\weight_i^2(\mbf{x})=1,\qquad 
\left\{\begin{array}{l}
\vspace{1mm}
\weight_i^1(\mbf{x})\ge \weight_i^2(\mbf{x})\quad\T{if }(B_i\mbf{n}|_{\partial\Omega^1})(\mbf{x})\ge 0,\\
\weight_i^1(\mbf{x})< \weight_i^2(\mbf{x})\quad\T{otherwise.}
\end{array}\right.
\end{equation}

We let $\pjump:\mathbb{R}^{2n}\rightarrow \mathbb{R}^n$ denote the function describing the diffusive flux across the interface $\Gatr$, that is
\begin{equation}\label{eq:pjump}
\pjump(\bx^1,\bx^2)
:=\mbf{P}(\bx^1,\bx^2)
(\bx^1-\bx^2)\quad\forall\, \bx^1,\bx^2\in\mathbb{R}^n,
\end{equation}
and assume that $\pjump\in C^{1,1}(\mathbb{R}^{2n})$ and that its Jacobian $\pjump'$ is bounded.

Throughout this work, we assume that the problem given by~\eqref{pde}, \eqref{modelbc}, \eqref{bcvec}, and \eqref{transmission} has a unique solution that remains bounded up to, and including, the final time $T$.

\section{Space discretisation by the discontinuous Galerkin method}\label{dg_heat}

\subsection{Finite element spaces}
Let $\mathcal{T}$ be a shape-regular and locally quasi-uniform subdivision of $\Omega$ into disjoint open elements $\el$, such that $\Gatr\subset \cup_{\el}\partial\kappa=:\Gamma$, the \emph{skeleton}. Further we decompose $\Gamma$ into three disjoint subsets $\Gamma=\partial\Omega\cup\Gamma_{\dint}\cup\Gatr$, where $\Gamma_{\dint}:=\Gamma\backslash(\partial\Omega\cup\Gatr)$. We assume that the subdivision $\mathcal{T}$ is constructed via mappings $F_{\kappa}$, where $F_{\kappa}:\hat{\kappa}\to\kappa$ are smooth maps with non-singular Jacobian, and $\hat{\kappa}$ is the reference $d$-dimensional simplex or the reference $d$-dimensional (hyper)cube. It is assumed that the union of the closures of the elements $\el$ forms a covering of the closure of $\Omega$; i.e., $\bar{\Omega}=\cup_{\el}\bar{\kappa}$. 
 
For $m \in \mathbb N$ we denote by $\mathbb{P}_m(\hat{\kappa})$ the set of polynomials of total degree at most $m$ if $\hat{\kappa}$ is the reference simplex, and the set of all tensor-product polynomials on $\hat{\kappa}$ of degree $k$ in each variable, if $\hat{\kappa}$ is the reference hypercube.  Let $m_{\kappa}\in\mathbb{N}$ be given for each $\el$.  We consider the $hp$-discontinuous finite element space
\begin{equation}\label{eq:FEM-spc}
\FEspace:=\{v\in L^2(\Omega):v|_{\kappa}\circ F_{\kappa}
   \in \mathbb{P}_{m_{\kappa}}(\hat{\kappa}),\,\el\},
\end{equation}
and set $\FESpace:=[\FEspace]^n$.

Next, we introduce relevant trace operators. Let $\kappa^+$, $\kappa^-$ be two elements sharing an edge $e:=\partial\kappa^+\cap\partial\kappa^-\subset\Gamma_{\dint}\cup \Gatr$. Denote the outward normal unit vectors on $e$ of $\partial\kappa^+$ and $\partial\kappa^-$ by $\mbf{n}^+$ and $\mbf{n}^-$, respectively. For functions $\mbf{q}:\Omega\to\mathbb{R}^n$ and $\mbf{Q}:\Omega\to\mathbb{R}^{n\times d}$ that may be discontinuous across $\Gamma$, we define the following quantities: for $\mbf{q}^+:=\mbf{q}|_{\kappa^+}$, $\mbf{q}^-:=\mbf{q}|_{\kappa^-}$ and $\mbf{Q}^+:=\mbf{Q}|_{\kappa^+}$, $\mbf{Q}^-:=\mbf{Q}|_{\kappa^-}$ on the restriction to $e$, we set
$ \mean{\mbf{q}}:=\ha(\mbf{q}^+ + \mbf{q}^-)$, $ \mean{\mbf{Q}}:=\ha(\mbf{Q}^+ + \mbf{Q}^-)$,
and
$ \jumptwo{\mbf{q}}:=\mbf{q}^+\otimes \mbf{n}^++\mbf{q}^-\otimes \mbf{n}^-$,
$ \jump{\mbf{Q}}:=\mbf{Q}^+ \mbf{n}^++\mbf{Q}^- \mbf{n}^-$,
where $\otimes$ denotes the standard tensor product operator, with $\mbf{q}\otimes\mbf{w}=\mbf{q}\mbf{w}^{{\mathsf T}}$. If $e\in \partial\kappa\cap\Gamma_{\partial}$, these definitions are modified as follows: $\mean{\mbf{q}}:=\mbf{q}^+ ,\ \mean{\mbf{Q}}:=\mbf{Q}^+$ and $\jumptwo{\mbf{q}}:=q^+\otimes\mbf{n},\ \jump{\mbf{Q}}:=\mbf{Q}^+\mbf{n}$.

Further, we introduce the mesh quantities $\tth:\Omega\to \mathbb{R}$, ${\tt m}:\Omega\to \mathbb{R}$ with $\tth(x)= \diam{\kappa}$, ${\tt m}(x)= m_{\kappa}$, if $x\in \kappa$, and the averaged values $\tth(x)= \mean{\tth}$, ${\tt m}(x)= \mean{{\tt m}}$, if $x\in \Gamma$.  Finally, we define $h_{\max}:=\max_{x\in{\Omega}}\tth$ and $h_{{\rm min}}:=\min_{x\in{\Omega}}\tth$.

We shall assume the existence of a constant $C_A\ge1$ independent of $\mathcal{T}$ such that, on any face  that is \emph{not} contained in $\Gatr$, given the two elements $\kappa$, $\kappa'$ sharing that face, the diffusion matrix $A$ satisfies
\begin{equation}\label{blv_A}
C_A^{-1} \le \big\| A\|_{\infty,\kappa}\|A^{-1}\big\|_{\infty,\kappa'}\le C_A.
\end{equation}
This assumption can be removed using the ideas from \cite{MR2491426}, but we refrain from doing so here for simplicity of the presentation.

The next result is a modification of the classical trace estimate for functions in $H^1(\Omega^1\cup\Omega^2)+\FEspace$; see~\cite{BuffaOrtner} for similar results.

\begin{lemma}\label{broken_trace}
Assume that the mesh $\mathcal{T}$ is both shape-regular and locally quasi-uniform. Then for $v\in H^1(\Omega^1\cup\Omega^2)+\FEspace$, the following trace estimate holds:
\begin{equation}\label{trace_estimate}
\sum_{j=1}^2\ltwo{v|_{\Omega^j}}{\Gatr}^2\le c_1 \epsilon \Big(\su\ltwo{\nabla v}{\kappa}^2 + \ltwo{\tth^{-1/2}\jump{v}}{\Gint}^2\Big)+ c_2 \epsilon^{-1}\ltwo{v}{}^2,
\end{equation}
for any $\epsilon>h_{\max}$ and for some constants $c_1>0$ and $c_2>0$, depending only on the shape-regularity of the mesh and on the domain $\Omega$.
\end{lemma}
\begin{proof}
We use the decomposition of $v\in H^1(\Omega^1\cup\Omega^2)+\FEspace$ into a conforming part $v_c$ and a non-conforming part $v_d:=v-v_c\in \FEspace$. This decomposition is described in \cite{karakashian_pascal,karakashian_pascal_convergence} for functions in $\FEspace$; the extension to $H^1(\Omega^1\cup\Omega^2)+\FEspace$ follows by taking $v_c\in H^1(\Omega^1\cup\Omega^2)$. Using Theorem 2.1(iii) from \cite{karakashian_pascal_convergence}, there exists a $v_c^i\in  H^1(\Omega^1\cup\Omega^2)$, $i=1,2$, such that 
\begin{equation}\label{recovery}
\su\ltwo{D^{\alpha}(v-v_c^i)}{\kappa\cap\Omega^i}^2
\le
C \ltwo{\tth^{1/2-|\alpha |}\jump{v}}{\Gint\cap\Omega^i}^2,
\end{equation}
where $D^{\alpha}$ is the differentiation operator for a multi-index $\alpha$ with $|\alpha|=0,1$. Hence, (\ref{recovery}) implies
\begin{equation}\label{recovery_transmission}
\su\ltwo{D^{\alpha}(v-v_c)}{\kappa}^2
\le
 C \ltwo{\tth^{1/2-|\alpha |}\jump{v}}{\Gint}^2,
\end{equation}
for $(v_c)|_{\Omega^j}:=v_c^j$, $j=1,2$.

The triangle inequality implies
\begin{equation}\label{triangle_in}
\sum_{j=1}^2\ltwo{v|_{\Omega^j}}{\Gatr}^2
\le
2\sum_{j=1}^2\big(\ltwo{v_c|_{\Omega^j}}{\Gatr}^2+\ltwo{(v-v_c)|_{\Omega^j}}{\Gatr}^2\big).
\end{equation}
To bound the first term on the right-hand side of (\ref{triangle_in}), we note that $v_c\in H^1(\Omega^1\cup\Omega^2)$, giving
\begin{equation}\label{conforming_part}
\begin{aligned}
\sum_{j=1}^2\ltwo{v_c|_{\Omega^j}}{\Gatr}^2
&\le
C \sum_{j=1}^2 \big( \ltwo{v_c}{\Omega^j}^2 + \ltwo{v_c}{\Omega^j}\ltwo{\nabla v_c}{\Omega^j}\big)^{1/2}
\le
C \big( \ltwo{v_c}{}^2 + \ltwo{v_c}{}\ltwo{\nabla v_c}{}\big)^{1/2}\\
&\le
C\big(  \epsilon^{-1}\ltwo{v_c}{}^2 + \epsilon\ltwo{\nabla v_c}{}^2\big)^{1/2},
\end{aligned}
\end{equation}
for any $\epsilon>0$ sufficiently small.  To further bound the right-hand side of (\ref{conforming_part}), we use the triangle inequality for each term, viz.,
\begin{equation}\label{triangle_two}
\ltwo{v_c}{}\le \ltwo{v-v_c}{} +\ltwo{v}{}
\quad\text{and}\quad
\ltwo{\nabla v_c}{}^2\le 2\su\big(\ltwo{\nabla(v-v_c)}{\kappa}^2 +\ltwo{\nabla v}{\kappa}^2\big),
\end{equation}
in conjunction with (\ref{recovery_transmission}), to arrive at
\begin{equation}\label{conforming_part_two}
\begin{aligned}
\sum_{j=1}^2\ltwo{v_c|_{\Omega^j}}{\Gatr}^2
&\le
C\Big(  \epsilon^{-1} \ltwo{\tth^{1/2}\jump{v}}{\Gint}^2 +  \epsilon^{-1} \ltwo{v}{}^2  +  \epsilon \ltwo{\tth^{-1/2}\jump{v}}{\Gint}^2 +\epsilon \su\ltwo{\nabla v}{\kappa}^2\Big)^{1/2}\\
&\le
C\Big(   \epsilon^{-1} \ltwo{v}{}^2  +  \epsilon
\ltwo{\tth^{-1/2}\jump{v}}{\Gint}^2 +\epsilon \su\ltwo{\nabla v}{\kappa}^2\Big)^{1/2},
\end{aligned}
\end{equation}
using the assumption that $\epsilon>h_{\max}$.

To bound the second term on the right-hand side of (\ref{triangle_in}), we use the trace estimate on each element in conjunction with (\ref{recovery_transmission}), viz.,
\begin{equation}\label{nonconforming_part}
\begin{aligned}
\sum_{j=1}^2\ltwo{(v-v_c)|_{\Omega^j}}{\Gatr}^2
&\le
C\!\!\!\!\!\!\!\!\sum_{\el: \bar{\kappa}\cap\Gatr\ne \emptyset}\!\!\!\!\!\!\!\!\big(     \ltwo{\tth^{-1/2}(v-v_c)}{\kappa}^2  +   \ltwo{\tth^{1/2}\nabla(v-v_c)}{\kappa}^2\big)
\le
C  \ltwo{\jump{v}}{\Gint}^2.
\end{aligned}
\end{equation}
Noting the assumption $\epsilon>h_{\max}$, the use of (\ref{conforming_part_two}) and (\ref{nonconforming_part}) in (\ref{triangle_in}) concludes the proof.
\end{proof}

\subsection{Space discretization}

The discretization of the space variables will be based on a discontinuous Galerkin method of interior penalty type for the diffusion part and of upwind type for the advection. Special care has to be given to the incorporation of the interface conditions. 

More specifically, we shall introduce a \emph{interior penalty discontinuous Galerkin} (IPDG, for short)  discretisation of the advection-diffusion operator 
\begin{equation}\label{eq:sec-system}
-\nabla\cdot (A \nabla\bw-W\mbf{B})
,
\end{equation}
where $W:=\diag(w_1,w_2,\dots, w_n)$, for $\bw:=(w_1,w_2,\dots, w_n)^{{\mathsf T}}$. 

To this end, we define $B(\bw_h,\bv_h)$ to be the IPDG bilinear form 
\begin{equation}\label{eq:bilinear}
\begin{aligned}
&\su\int_{\kappa} (A\nabla\bw_h-W_h\mbf{B}): \nabla \bv_h
+\int_{\Gatr}\Big(
\mean{W_h\mbf{B}}+\mathcal{B}_{\inter}\jumptwo{\bw_h}
\Big):\jumptwo{\bv_h}\\
&-\int_{\Gint}\Big(\mean{A\nabla \bw_h-W_h\mbf{B}}:\jumptwo{\bv_h}+\mean{A\nabla \bv_h}:\jumptwo{\bw_h}
-(\Sigma+\mathcal{B})\jumptwo{\bw_h}:\jumptwo{\bv_h}\Big)\\
&-\int_{\Gamma_{\ddd}}\Big((A\nabla \bw_h-\mathcal{X}^+W_h\mbf{B}):(\bv_h\otimes \mbf{n})+(A\nabla \bv_h):(\bw_h\otimes \mbf{n})
-\Sigma\bw_h\cdot\bv_h\Big)\\
&+\int_{\Gamma_{\dn}}(\mathcal{X}^+W_h\mbf{B}):(\bv_h\otimes \mbf{n}),
\end{aligned}
\end{equation}
and define
\begin{equation}\label{eq:nonlinear}
N(\bw_h,\bv_h):=
\int_{\Gatr}\left(\mbf{P}(\bw_h)\jumptwo{\bw_h}
-({\rm I}-\Friction)
\left(
\mean{W_h\mbf{B}}+\mathcal{B}_{\inter}\jumptwo{\bw_h}
\right)
\right):\jumptwo{\bv_h}.
\end{equation}
Here, $\Sigma:=C_{\sigma}A{\tt m}^2\tth^{-1}$ denotes the \emph{discontinuity-penalization parameter} matrix with $C_{\sigma}>1$ constant.  Furthermore, $\mathcal{B}:=\frac{1}{2}\diag(|B_1\cdot \mbf{n}|,\dots,|B_n\cdot \mbf{n}|)$ and $\mathcal{B}_{\inter}:=(\Weight^1-\frac{1}{2}{\rm I})\mbf{B}\mbf{n}|_{\Omega^1}= (\Weight^2-\frac{1}{2}{\rm I})\mbf{B}\mbf{n}|_{\Omega^2}$ is diagonal with non negative entries.

\begin{remark} We comment on the interface terms. The diffusion term, appearing in $N$, is simply given by $\int_{\Gatr}\mbf{P}(\bw)\jumptwo{\bw}:\jumptwo{\bv_h}$. This resembles the typical jump stabilisation term with the permeability coefficient replacing the discontinuity-penalization parameter, rendering its implementation within a dG computer code straightforward.

To ensure the coercivity of $B$, the advective interface term has been split as $\Friction ={\rm I}-({\rm I}-\Friction)$, resulting into contributions in both $B$ and $N$. Indeed, the advective interface contribution in $B$ can be recast using the weighted mean $\mean{W_h\mbf{B}}^\weight:=\Weight^1W_h\mbf{B}|_{\Omega^1} + \Weight^2W_h\mbf{B}|_{\Omega^2}$, so that
\begin{equation}\label{eq:weighted_average}
\mean{W_h\mbf{B}}^\weight:\jumptwo{\bv_h}=\left(
\mean{W_h\mbf{B}}+\mathcal{B}_{\inter}\jumptwo{\bw_h}
\right):\jumptwo{\bv_h},
\end{equation}
thereby resembling the typical dG upwinding for linear advection problems and, hence, ensuring the coercivity of $B$.
\end{remark}
\begin{remark}
In this setting, $\Gamma_{\dn}$ can have non-trivial intersection with both $\partial\Omega^-_i$ and $\partial\Omega^+_i$, thereby extending the discontinuous Galerkin method proposed in \cite{houston-etal:01}.
\end{remark}
  
\section{Elliptic projection error}\label{err_ell_problem}
The a priori error analysis is based on a (non-standard) elliptic projection inspired by a construction of Douglas and Dupont for the treatment of nonlinear boundary conditions in the context of conforming finite element methods \cite{douglas_dupont}. 

\begin{definition}
For each $t\in [0,T]$ we define the \emph{elliptic projection} $\bw_h\in \FESpace$ to be the solution of the problem: find $\bw_h\equiv\bw_h(t)\in \FESpace$, such that 
\begin{equation}\label{ellipitc_projection}
{B}(\bu-\bw_h,\bv_h)+\lambda \ltwoin{\bu-\bw_h}{\bv_h}
+N(\bu,\bv_h)-N(\bw_h,\bv_h)
=
0\quad \forall \bv_h\in \FESpace,
\end{equation}
for some fixed $\lambda>0$, and $\bu$ denoting the exact solution. 
\end{definition}

The constant $\lambda>0$ in the definition above will be chosen large enough to ensure the uniqueness of the projection $\bw_h$, cf. Lemma \ref{ell_proj_lemma} below.

Next, denoting by $\SSpace^s:=\Hn^s+\FESpace$, $s\in\mathbb{R}$,
we define the dG-norm on $\SSpace^1$
\begin{equation}
\label{eq:dGnorm}
\begin{aligned}
\ndg{\bw}:=\bigg(&\su\big(\ltwo{\sqrt{A}\nabla\bw}{\kappa}^2+\frac{1}{2}\ltwo{\sqrt{\diag(\nabla\cdot \mbf{B})}\bw}{\kappa}^2\big)
+\ltwo{\sqrt{\Sigma}\jumptwo{\bw}}{\Gamma_{\ddd}\cup\Gamma_{\dint}}^2
\\
&+\ltwo{\sqrt{\mathcal{B}}\jumptwo{\bw}}{\Gamma\backslash\Gatr}^2
+\ltwo{
\sqrt{\mathcal{B}_{\inter}}
\jumptwo{\bw}}{\Gatr}^2
\bigg)^{1/2},
\end{aligned}
\end{equation}
where $
\ltwo{\mbf{Q}}{\kappa}^2 := \int_\kappa \sum_{i=1}^n|Q_i(x)|^2 \, {\rm d}x$,
denotes the Frobenius norm whenever $\mbf{Q}$ is a $n\times d$ tensor. We assume that~\eqref{eq:dGnorm} is a norm. This is satisfied when standard assumptions on the solution in conjunction with the boundary conditions hold on each subdomain, e.g., $\Gamma_{\ddd}\cap\partial\Omega^j$ has positive $(d-1)$-dimensional (Hausdorff) measure for $j=1,2$. If the interface manifold $\Gatr$ is not characteristic to the advection field, such hypotheses can be further relaxed.

For the remainder of this work, we shall make the simplifying assumption that $\mbf{B}$ is such that:
\begin{equation}\label{wind_assumption}
 B_i\cdot\nabla (v_h)_i\in \FEspace,\quad \text{ for }\ i=1,\dots,n,
\end{equation}
for any function $\bv_h:=((v_h)_1,\dots,(v_h)_n)^{{\mathsf T}}\in \FESpace$. We note, however, that this appears not to be a genuine limitation in the arguments presented below: ideas on how to circumvent this assumption have been presented, e.g., in \cite{houston-etal:01, ayuso_marini}, for the case of scalar linear advection-diffusion problems. 

The next two results show the coercivity and the continuity of the bilinear form $B(\cdot,\cdot)$. Their proofs follow straightforward variations of well-known arguments (see, e.g., \cite{arnold,houston-etal:01}) and are, therefore, omitted for brevity.

\begin{lemma}\label{coercivity}
There exists a positive constant $C_{coer}\in\mathbb{R}$, such that, for $\bv_h\in \FESpace$,
\[
B(\bv_h,\bv_h)\ge C_{coer}\ndg{\bv_h}^2.
\]
\end{lemma}
\begin{lemma}\label{continuity}
Let $\Pi:[L^2(\Omega)]^n\to \FESpace$ denote the $L^2$-orthogonal projection onto $\FESpace$.
For any $\bw\in \Hn^{s}$, $s>3/2$ and $\bv_h\in \FESpace$ we have
\[
|{B}(\boeta,\bv_h)|
\le   C_{cont} \ndg{\boeta}_{\mathcal{B}}
\ndg{\bv_h},
\]
where $\boeta:=\bw-\Pi \bw$ and 
\begin{equation}\label{eq:Aeta}
\ndg{\boeta}^2_{\mathcal{B}}:=
\ndg{\boeta}^2
+\ltwo{{\Sigma}^{-1/2}\mean{A\nabla \boeta}}{\Gamma_{\ddd}\cup\Gamma_{\dint}}^2
+\ltwo{\sqrt{\mathcal{B}}\mean{\boeta}}{\Gamma}^2
+\ltwo{\sqrt{\mathcal{B}_{\inter}}\mean{\boeta}}{\Gatr}^2.
\end{equation}
\end{lemma}

The next result establishes the well-posedness of the
problem~\eqref{ellipitc_projection} and relevant approximation properties.

\begin{lemma}\label{ell_proj_lemma}
Assume that $\bu\in \Hn^s$, $s>3/2$
for all $t\in(0,T]$.
For $\lambda>0$ sufficiently large and for $h_{\max}$ sufficiently small, the variational problem (\ref{ellipitc_projection}) has a unique solution $\bw_h\in \FESpace$ for each $t\in(0,T]$. Moreover, the following bound holds:
\begin{equation}
\label{ell_proj_bound}
C_{coer}\ndg{\brho}^2+ \lambda \ltwo{\brho}{}^2
\le \ndg{\boeta}_{\mathcal{B},\lambda}^2,
\end{equation}
and, if also $\bu_t\in \Hn^s$, then
\begin{equation}\label{ell_proj_bound_t}
C_{coer}\ndg{\brho_t}^2+ \lambda \ltwo{\brho_t}{}^2 \le \ndg{\boeta_t}_{\mathcal{B},\lambda}^2+\ndg{\boeta}_{\mathcal{B},\lambda}^2
\end{equation}
where $\brho:=\bu-\bw_h$, $\boeta:=\bu-\Pi\bu$, and
\[
\ndg{\boeta}_{\mathcal{B},\lambda}^2:=
C_{c}\ndg{\boeta}^2_{\mathcal{B}}+7\lambda \ltwo{\boeta}{}^2,
\]
with $C_{c}:=(4C_{cont}^2+3C_{coer}^2)/C_{coer}$.
\end{lemma}

{\bf Proof.}
Well-posedness of (\ref{ellipitc_projection}) is established by proving that the associated mapping is strongly monotone on $\FESpace$.  Using the assumption that $\pjump\in C^{0,1}(\mathbb{R}^{2n})$
, we  get 
\begin{equation}\label{eq:N_bound}
\begin{aligned}
|N(\bv,\bz)-N(\bw,\bz)|&\le \int_{\Gatr} (|\pjump (\bv)-\pjump (\bw)|+C_\mathcal{B} |\bv-\bw)|)\,|\jumptwo{\bz}|\\
&\le \int_{\Gatr} \left(C^p|\bv-\bw|+C_\mathcal{B}|\bv-\bw|\right) \,|\jumptwo\bz|\\
&\le C^{p}_{\mathcal{B}} \sum_{j=1}^2(\ltwo{(\bv-\bw)|_{\Omega^j}}{\Gatr}^2 +\ltwo{\bz|_{\Omega^j}}{\Gatr}^2),
\end{aligned}
\end{equation}
where $C^p$ is a Lipschitz constant for the function $\pjump$,  $C_{\mathcal{B}}>0$ is a  constant proportional to $\displaystyle{\max_{i=1,\dots,n}\ltwo{B_i\mbf{n}}{\infty,\Gatr}}$, and $C^{p}_{\mathcal{B}}>0$ is a constant depending on both $C^p$ and $C_{\mathcal{B}}$.

This, in conjunction with the coercivity of the bilinear form $B$ provided by Lemma~\ref{coercivity}, gives
\[
\begin{aligned}
&B(\bv_h-\bw_h,\bv_h-\bw_h)+\lambda \ltwoin{\bv_h-\bw_h}{\bv_h-\bw_h}
+N(\bv_h,\bv_h-\bw_h)-N(\bw_h,\bv_h-\bw_h)\\
&\ge C_{coer}\ndg{\bv_h-\bw_h}^2
+\lambda \ltwo{\bv_h-\bw_h}{}^2-C^p_{\mathcal{B}}\sum_{j=1}^2\ltwo{(\bv_h-\bw_h)|_{\Omega^j}}{\Gatr}^2\\
&\ge \frac{C_{coer}}{2}\ndg{\bv_h-\bw_h}^2
+\big(\lambda
-\frac{2 C^p_\mathcal{B} c_1c_2}{C_{coer}\alpha_{\min}}\big)\ltwo{\bv_h-\bw_h}{}^2,
\end{aligned}
\]
the last inequality owning to the trace estimate~\eqref{trace_estimate} with $\epsilon= C_{coer}\alpha_{\min}/(2C^p_{\mathcal{B}}c_1)$. Hence strong monotonicity is ensured as soon as $\lambda>2 C^p_\mathcal{B}c_1c_2/C_{coer}\alpha_{\min}$.

To show~\eqref{ell_proj_bound} and~\eqref{ell_proj_bound_t}, we split the error $\bu-\bw_h=\boeta+\bxi$, with $\boeta:=\bu-\Pi\bu$ and $\bxi:=\Pi\bu-\bw_h$.  Then, setting $\bv_h=\bxi$ in (\ref{ellipitc_projection}), we deduce
\begin{equation}\label{error_eq_ell_proj}
B(\bxi,\bxi)+\lambda \ltwoin{\bxi}{\bxi}= -{B}(\boeta,\bxi)-\lambda \ltwoin{\boeta}{\bxi}-N(\bu,\bxi)+N(\bw_h,\bxi).
\end{equation}
Using the coercivity and the continuity of $B$ along with the bound~\eqref{eq:N_bound} in (\ref{error_eq_ell_proj}) gives
\begin{equation}\label{error_in_elliptic}
 \frac{3}{4}C_{coer}\ndg{\bxi}^2+ \frac{3}{4}\lambda \ltwo{\bxi}{}^2\le\frac{C_{cont}^2}{C_{coer}}\ndg{\boeta}^2_{\mathcal{B}}+\lambda \ltwo{\boeta}{}^2+C^p_\mathcal{B}\sum_{j=1}^2\big(\ltwo{\boeta|_{\Omega^j}}{\Gatr}^2+\ltwo{\bxi|_{\Omega^j}}{\Gatr}^2\big).
\end{equation}
Using (\ref{trace_estimate}) with
$\epsilon= C_{coer}\alpha_{\min}/(4C^p_\mathcal{B}c_1)$ on the last term on the right-hand side of (\ref{error_in_elliptic}), we arrive at
\[
 \frac{1}{2}C_{coer}\ndg{\bxi}^2+ \frac{3}{4}\lambda \ltwo{\bxi}{}^2\le 
C_{c}\ndg{\boeta}^2_{\mathcal{B}}+\lambda \ltwo{\boeta}{}^2
  +C^p_\mathcal{B}c_2 \epsilon^{-1}(\ltwo{\boeta}{}^2
+\ltwo{\bxi}{}^2).
\]
Choosing $\lambda> 16(C^p_\mathcal{B})^2c_1c_2/(C_{coer}\alpha_{\min})$, we deduce
\begin{equation}\label{error_in_elliptic_three}
 \frac{1}{2}C_{coer}\ndg{\bxi}^2+ \frac{1}{2}\lambda \ltwo{\bxi}{}^2\le 
C_{c}\ndg{\boeta}^2_{\mathcal{B}}+\frac{5}{4}\lambda \ltwo{\boeta}{}^2.
\end{equation}
A triangle inequality already gives~\eqref{ell_proj_bound}.

In view of obtaining~\eqref{ell_proj_bound_t}, we differentiate~\eqref{ellipitc_projection} with respect to $t$ and then test with $\bv_h=\bxi_t$:
\begin{equation}\label{error_in_elliptic_four}
\begin{aligned}[l]
&B\left(\bxi_t,\bxi_t\right)+\lambda \ltwoin{\bxi_t}{\bxi_t}= 
-{B}\left(\boeta_t,\bxi_t\right)-\lambda \ltwoin{\boeta_t}{\bxi_t}\\
&\quad -\int_{\Gatr}
\frac{d}{d t}
\big(
\pjump (\bu)-\pjump (\bw_h)
+({\rm I}-\Friction)
\left(
\mean{(U-W_h)\mbf{B}}+\mathcal{B}_{\inter}\jumptwo{\brho}
\right)
\big):\jumptwo{\bxi_t}.
\end{aligned}
\end{equation}
Using the assumption that the Jacobian $\pjump'\in C^{0,1}(\mathbb{R}^{2n})$ and is bounded we obtain
\begin{equation}\label{eq:p_t_bound}
\begin{aligned}
\left|\int_{\Gatr}\right.&\left.  \frac{d}{d t} \big( \pjump (\bu)-\pjump (\bw_h) +({\rm I}-\Friction) \left( \mean{(U-W_h)\mbf{B}}+\mathcal{B}_{\inter}\jumptwo{\brho} \right) \big):\jumptwo{\bxi_t}\right|\\
&\le \int_{\Gatr}\Big(|(\pjump'(\bu)-\pjump'(\bw_h))[\bu_t]| +|\pjump'(\bw_h)[\brho_t]|+C_{\mathcal{B}}\sum_{j=1}^2(|\brho_t|_{\Omega^j}|+|\brho|_{\Omega^j}|)\Big) |\jumptwo{\bxi_t}| \\
&\le \Big(\sum_{j=1}^2(\tilde{C}^{p}_{\mathcal{B},\bu} \ltwo{\brho|_{\Omega^j}}{\Gatr} +C^{p}_{\mathcal{B}}     \ltwo{\brho_t|_{\Omega^j}}{\Gatr})\Big) \sum_{j=1}^2\ltwo{\bxi_t|_{\Omega^j}}{\Gatr}\\
&\le \sum_{j=1}^2\left(\tilde{C}^{p}_{\mathcal{B},\bu} \ltwo{\brho|_{\Omega^j}}{\Gatr}^2 +2C^{p}_{\mathcal{B}} \ltwo{\boeta_t|_{\Omega^j}}{\Gatr}^2 +(\tilde{C}^{p}_{\mathcal{B},\bu} +3C^{p}_{\mathcal{B}})\ltwo{\bxi_t|_{\Omega^j}}{\Gatr}^2\right) \\
&\le  C^{p}_{\mathcal{B},\bu}\sum_{j=1}^2\left( \ltwo{\brho|_{\Omega^j}}{\Gatr}^2 +\ltwo{\boeta_t|_{\Omega^j}}{\Gatr}^2 +\ltwo{\bxi_t|_{\Omega^j}}{\Gatr}^2\right), 
\end{aligned}
\end{equation}
with $\tilde{C}^{p}_{\mathcal{B},\bu}$, $C^{p}_{\mathcal{B},\bu}>0$ constants depending on $C^{p}_{\mathcal{B}}$ and on $\ltwo{\bu_t}{L^{\infty}([0,T]\times \Omega)}$.  Here and in what follows, square brackets are used to denote the argument of a linear operator.

Applying~\eqref{eq:p_t_bound} to the right-hand side of~\eqref{error_in_elliptic_four} and using the coercivity and continuity of $B$ yields
\begin{equation}\label{error_in_elliptic_five}
\begin{aligned}
 \frac{3}{4}C_{coer}\ndg{\bxi_t}^2+ \frac{3}{4}\lambda \ltwo{\bxi_t}{}^2&\le 
\frac{C_{cont}^2}{C_{coer}} \ndg{\boeta_t}^2_{\mathcal{B}}
+\lambda \ltwo{\boeta_t}{}^2\\
&+C^{p}_{\mathcal{B},\bu}\sum_{j=1}^2\left(
\ltwo{\brho|_{\Omega^j}}{\Gatr}^2
+\ltwo{\boeta_t|_{\Omega^j}}{\Gatr}^2
+\ltwo{\bxi_t|_{\Omega^j}}{\Gatr}^2\right).
\end{aligned}
\end{equation}

As before, using the trace estimate (\ref{trace_estimate}) on the last term on the right-hand side of (\ref{error_in_elliptic_five}) with
$\epsilon= C_{coer}\alpha_{\min}/(4C^{p}_{\mathcal{B},\bu}c_1)$, we arrive to
\[
\begin{aligned}
 \frac{1}{2}C_{coer}\ndg{\bxi_t}^2+ \frac{1}{2}\lambda \ltwo{\bxi_t}{}^2&\le 
\frac{4C_{cont}^2+3C_{coer}^2}{C_{coer}} \ndg{\boeta_t}^2_{\mathcal{B}}
        +\frac{5}{4}\lambda \ltwo{\boeta_t}{}^2
+\frac{1}{4}\ndg{\boeta}_{\mathcal{B},\lambda}^2,
\end{aligned}
\]
for any $\lambda> 16(\max\{C^{p}_{\mathcal{B},\bu},C^{p}_{\mathcal{B}}\})^2c_1c_2/(C_{coer}\alpha_{\min})$.  Now~\eqref{ell_proj_bound_t} easily follows by the triangular inequality.
\qed

We conclude this section with an $L^2$-error bound of the elliptic projection~(\ref{ellipitc_projection}). This is obtained by an Aubin-Nitsche duality-type argument, inspired by a construction of Douglas and Dupont~\cite{douglas_dupont} for nonlinear boundary conditions.

The interface operator $N$ given in~\eqref{eq:nonlinear} consists of a nonlinear component driven by the function $\pjump(\bw)=\mbf{P}(\bw)\jumptwo{\bw}$ and a linear component which we can characterise by  introducing the linear operator $L[\bw]:=-({\rm I}-\Friction) \mean{W\mbf{B}}^\weight +\mathcal{B}_{\inter}\jumptwo{\bw}$. We abbreviate $\SSpace := \SSpace^1$ and let $\SSpace^*$ be the dual space of $\SSpace$. It is convenient to momentarily view $N$ as an operator from $\SSpace \to \SSpace^*$, indicated with a calligraphic font:
\begin{equation*}
\mathcal N : \; \SSpace \to \SSpace^*, \; \bw \mapsto \Bigl( \bv \mapsto \int_{\Gatr}\left(\pjump(\bw) 
+L[\bw]
\right):\jumptwo{\bv} \Bigr),
\end{equation*}
where the dependence on $\bv$ represents a linear mapping $\SSpace \to \mathbb R$ in $\SSpace^*$. Thus the derivative $\mathcal N'$ is a mapping $\SSpace \to L(\SSpace,\SSpace^*)$, where $L(\SSpace,\SSpace^*)$ denotes the linear mappings from $\SSpace$ to $\SSpace^*$. Therefore the integral
\[
\mathcal{P}(t,\bv) := \int_0^1 \mathcal N'(\bw^\theta(t,\cdot))(\bv)\, d\theta,
\]
where $\bw^\theta := \theta \bu + (1-\theta) \bw_h$,
belongs to $\SSpace^*$ for each $t \in (0,T)$, $\bv \in \SSpace$. In particular $\mathcal{P}(t, \bu(t,\cdot)-\bw_h(t,\cdot)) \in \SSpace^*$ and
\begin{align}\nonumber
\mathcal{P}(t, \bu(t,\cdot)-\bw_h(t,\cdot)) \; 
= & \int_0^1 \mathcal N'(\bw^\theta(t,\cdot))(\bu(t,\cdot)-\bw_h(t,\cdot)) \, d\theta\\ 
= & \int_0^1 \partial_\theta (\mathcal N(\bw^\theta(t,\cdot))) \, d\theta= \, \mathcal N(\bu(t,\cdot)) - \mathcal N(\bw_h(t,\cdot)),\label{eq:PN}
\end{align}
using that $[0,1] \to \SSpace^*, \; \theta \mapsto \mathcal N(\bw^\theta(t,\cdot))$ is continuously differentiable as $\pjump \in C^{1,1}(\mathbb{R}^{2n})$. We shall frequently abbreviate $\mathcal{P}(t,\bz(t,\cdot))$ by $\mathcal{P}\bz$ below.

We assume that there is an $s \in (3/2,2]$ such that for all $\mbf{\alpha} \in [L^2(\Omega)]^n$ and $\mbf{\beta} \in [H^{1/2}(\Gatr)]^{2n}$ there exists a solution $\mbf{\zeta} \in \Hn^s$ of the linear dual equation:
\begin{equation}\label{eq:weak_dual_prob}
B(\bv,\mbf{\zeta})+\lambda \ltwoin{\bv}{\mbf{\zeta}}
+\langle{\mathcal{P}}\bv,{\mbf{\zeta}}\rangle
=(\bv,\mbf{\alpha}) + \langle\bv, \mbf{\beta}\rangle_{\Gatr} \quad \forall \, \bv\in\honen.
\end{equation}
Further, we assume that the dual solution $\mbf{\zeta}$ satisfies the elliptic regularity bound:
\begin{equation}\label{eq:dual_elliptic_regularity}
  \sum_{j=1}^2\ltwo{\mbf{\zeta}}{H^s(\Omega^j)}\lesssim\ltwo{\mbf \alpha}{} + \ltwo{\mbf \beta}{H^{1/2}(\Gatr)}.
\end{equation}
The bound is motivated by \cite{mitreamitrea} and \cite{savare}: Suppose $\Omega^1$ and $\Omega^2$ are smooth or creased domains. Assuming the existence of $\mbf{\zeta}$ allows to decouple the problem into the subdomains, using $\mbf{\zeta}$ for the boundary data on the interface. With $\honen$ control available, lower-order terms may be moved to the right-hand side before applying elliptic regularity bounds such as those mentioned above, leading typically to control in fractional Sobolev norms.

\begin{lemma}\label{ell_proj_l2_lemma}
Assume that the hypothesis of Lemma \ref{ell_proj_lemma} and~\eqref{eq:weak_dual_prob} with~\eqref{eq:dual_elliptic_regularity} hold true.
For $\lambda>0$ sufficiently large, for $h_{\max}$ sufficiently small, the following error bounds holds:
\begin{align}
\label{ell_proj_l2_bound}
\ltwo{\brho}{}
&\le  C  (1+h^2_{\max}\lambda)^{1/2} h^{s-1}_{\max}\,
	{\ndg{\boeta}}_{\mathcal{B},\lambda},
\end{align}
If in addition the function  
$\pjump$ 
defined in \eqref{eq:pjump} is twice differentiable with bounded second partials and $\bu, \bu_t \in W^{1,\infty}([0,T] \times \Gatr)$ then
\begin{align}
\label{ell_proj_l2_bound_t}
\ltwo{\brho_t}{}^2
&\le  C  (1+h^2_{\max}\lambda)^{1/2} h^{s-1}_{\max}\,({\ndg{\boeta_t}}_{\mathcal{B},\lambda}+{\ndg{\boeta}}_{\mathcal{B},\lambda}).
\end{align}
The constant $C$ depends only on $C_A$ and the shape-regularity of the mesh.
\end{lemma}

{\bf Proof.} Let $\bz$ solve \eqref{eq:weak_dual_prob} with $\mbf{\alpha} = \brho$ and $\mbf{\beta} = \mbf{0}$.
Owing to the adjoint consistency of the symmetric interior penality bilinear form $B$ as well as of $\mathcal P$, the
dual solution $\bz$ also satisfies \eqref{eq:weak_dual_prob} for $\bv \in \FESpace$: Because $\bz$ belongs to $\honen$, the integral terms over $\Gint$ in $B$ vanish, so that upon integration by parts $B(\bv,\bz)$ is equal to the $L^2(\Omega)$ scalar product of $\bv$ and the adjoint differential operator applied to $\bz$, as well as $L^2$ scalar products over $\Gatr$ and $\partial \Omega$. In all scalar products, also including those arising from $\mathcal P$, derivatives are then only acting on $\bz$ and not on $\bv$, implying by density of smooth functions in $L^2$ that \eqref{eq:weak_dual_prob} also holds for dG test functions $\bv$. This allows us to test in~\eqref{eq:weak_dual_prob} with $\bv=\brho=\bu-\bw_h$ and get:
\begin{equation}\label{eq:dual_l2}
\begin{aligned}
\ltwo{\brho}{}^2&={B}(\brho,\bz)+\lambda \ltwoin{\brho}{\bz}
+\langle{\mathcal{P}}\brho,\bz\rangle\\
&={B}(\brho,\bz)+\lambda \ltwoin{\brho}{\bz}
+N(\bu,\bz)-N(\bw_h,\bz)\\
&=
B(\brho,\boeta^z)
+\lambda \ltwoin{\brho}{\boeta^z}
+N(\bu,\boeta^z)-N(\bw_h,\boeta^z),
\end{aligned}
\end{equation}
with $\boeta^z=\bz-\Pi\bz$.  The second equality follows from \eqref{eq:PN}; the last equality  follows from the definition of the elliptic projection~(\ref{ellipitc_projection}).

We now bound each term on the right-hand side of~\eqref{eq:dual_l2}.
Observe that
\begin{equation}\label{eq:dual_B_bound}
\begin{aligned}
\left|
B(\brho,\boeta^z)\right|&\le C \ndg{\brho}(
\ndg{\boeta^z}_{\mathcal{B}}^2
+
\su (\max_{i=1,\dots,n}\ltwo{a_i^{-1/2}B_i}{\infty,\kappa}^2)\ltwo{\boeta^z}{\kappa}^2)^{1/2}\\
\vspace{2mm}
&\le C \ndg{\brho}(
\ndg{\boeta^z}_{\mathcal{B}}^2
+\lambda\ltwo{\boeta^z}{}^2)^{1/2}
\le C \ndg{\brho} \,{\ndg{\boeta^z}}_{\mathcal{B},\lambda},
\end{aligned}
\end{equation}
for $\lambda$ big enough.
The nonlinear interface term in~(\ref{eq:dual_l2}) is bounded
as in the proof of Lemma~\ref{ell_proj_lemma}, yielding  
\begin{equation}\label{eq:dual_inter_bound}
\begin{aligned}
&\left| N(\bu,\boeta^z)\right. \left. \hspace{-1mm}-N(\bw_h,\boeta^z)\right|
\le 4C^p_{\mathcal{B}}\big(\sum_{j=1}^2\ltwo{\brho|_{\Omega^j}}{\Gatr}^2\big)^{1/2}
\big(\sum_{j=1}^2\ltwo{\boeta^z|_{\Omega^j}}{\Gatr}^2\big)^{1/2}\\
\vspace{2mm}
&\quad\le \big(C_{coer}\ndg{\brho}^2+4C^p_{\mathcal{B}}c_2\epsilon^{-1}\ltwo{\brho}{}^2\big)^{1/2}
\big(C_{coer}\ndg{\boeta^z}^2+4C^p_{\mathcal{B}}c_2\epsilon^{-1}\ltwo{\boeta^z}{}^2\big)^{1/2}\\
&\quad\le \big(C_{coer}\ndg{\brho}^2+\lambda\ltwo{\brho}{}^2\big)^{1/2}
\big(C_{coer}\ndg{\boeta^z}^2+\lambda\ltwo{\boeta^z}{}^2\big)^{1/2},
\end{aligned}
\end{equation}
using once again (\ref{trace_estimate}) with
$\epsilon= C_{coer}\alpha_{\min}/(4C^p_\mathcal{B}c_1)$, 
for $\lambda> 16(C^p_\mathcal{B})^2c_1c_2/(C_{coer}\alpha_{\min})$.

Using the above inequalities 
and applying the bounds on the elliptic projection error given in
Lemma~\ref{ell_proj_lemma}, we obtain
\begin{equation}\label{rho_l2}
\ltwo{\brho}{}^2\le C
{\ndg{\boeta}}_{\mathcal{B},\lambda}
{\ndg{\boeta^z}}_{\mathcal{B},\lambda}
\end{equation}
for $\lambda$ big enough.

The term in $\boeta^z$ can be bounded using standard approximation estimates for the error of the (orthogonal) $L^2$-projection (see, e.g., \cite{schwab}) yielding
\begin{equation}\label{eta_z}
{\ndg{\boeta^z}}_{\mathcal{B},\lambda}^2
\le C\su h_\kappa^{2(s-1)}(1+h_\kappa^2\lambda) \, \ltwo{\bz}{H^s(\kappa)}^2.
\end{equation}

Further, inserting~(\ref{eta_z}) into~(\ref{rho_l2}) and then using \eqref{eq:dual_elliptic_regularity}
we get
\begin{equation}\label{rho_l2_1}
\begin{aligned}
	\ltwo{\brho}{}^2 &\le C(1+h^2_{\max}\lambda)^{1/2} h^{s-1}_{\max}\,
	{\ndg{\boeta}}_{\mathcal{B},\lambda}\sum_{j=1}^2\ltwo{\bz}{H^s(\Omega^j)}\\
	&\le C  (1+h^2_{\max}\lambda)^{1/2} h^{s-1}_{\max}\,
	{\ndg{\boeta}}_{\mathcal{B},\lambda}
	\ltwo{\brho}{},\end{aligned}    
\end{equation}
thus yielding~\eqref{ell_proj_l2_bound}.

We now prove~\eqref{ell_proj_l2_bound_t}. Let $\tilde{\bz}$ be the solution \eqref{eq:weak_dual_prob} with $\mbf{\alpha} = \brho_t$ and $\mbf{\beta} = \mbf{0}$.
Using the linearity of $\mathcal{P}$ in the second argument, we have
\begin{align*}
\textstyle d_t \mathcal{P}(t,\brho(t)) = \, & ( \partial_{\sbrho} \mathcal{P} ) \bigr|_{(t,\sbrho(t))} \brho_t + (\partial_t \mathcal{P})\bigr|_{(t,\sbrho(t))}
= \mathcal{P} \bigr|_{(t,\sbrho(t))} \brho_t + (\partial_t \mathcal{P})\bigr|_{(t,\sbrho(t))}.
\end{align*}
Testing in~\eqref{eq:weak_dual_prob} with $\bv=\brho_t$ gives, for each $t$,
\begin{equation}\label{eq:dual_l2_t}
\begin{aligned}
\ltwo{\brho_t(t)}{}^2&={B}(\brho_t,\tilde{\bz})+\lambda \ltwoin{\brho_t}{\tilde{\bz}}
+\langle{\mathcal{P}}\brho_t,{\tilde{\bz}}\rangle\\
&=
B(\brho_t,\tilde{\bz})+\lambda \ltwoin{\brho_t}{\tilde{\bz}}
+ \langle d_t  {\mathcal{P}}(t,\brho(t)),{\tilde{\bz}}\rangle- \langle \partial_t \mathcal{P}, {\tilde{\bz}} \rangle\\
&=
B(\brho_t,\tilde{\bz})+\lambda \ltwoin{\brho_t}{\tilde{\bz}}
+ \bigl\langle d_t \bigl( N(\bw,\cdot)-N(\bw_h,\cdot) \bigr), \tilde{\bz} \bigl\rangle  - \langle \partial_t \mathcal{P}, {\tilde{\bz}} \rangle\\
&=
B(\brho_t,\boeta^{\tilde{z}})+\lambda \ltwoin{\brho_t}{\boeta^{\tilde{z}}}
+ \bigl\langle d_t \bigl( N(\bw,\cdot)-N(\bw_h,\cdot) \bigr), \boeta^{\tilde{z}} \bigl\rangle - \langle \partial_t \mathcal{P}, {\tilde{\bz}} \rangle.
\end{aligned}
\end{equation}
with $\boeta^{\tilde{z}}=\tilde{\bz}-\Pi\tilde{\bz}$.  The last equality follows from differentiating with respect to time the definition of the elliptic projection~(\ref{ellipitc_projection}).

The first three terms on the right-hand side of~\eqref{eq:dual_l2_t} can be bounded by the same argument used above. In particular, following the argument in~\eqref{eq:p_t_bound} and then applying the bounds on the elliptic projection error given in Lemma~\ref{ell_proj_lemma} yields
\begin{equation}\label{eq:Ntbound}
\begin{aligned}
& |d_t \left(N(\bu,\boeta^{\tilde{z}})-N(\bw_h,\boeta^{\tilde{z}})\right)|\\
\le \; &
 4C^p_{\mathcal{B},\bu}\big(\sum_{j=1}^2(\ltwo{\brho|_{\Omega^j}}{\Gatr}^2+\ltwo{\brho_t|_{\Omega^j}}{\Gatr}^2)\big)^{1/2}
\big(\sum_{j=1}^2\ltwo{\boeta^{\tilde{z}}|_{\Omega^j}}{\Gatr}^2\big)^{1/2}\\
\le \; & C({\ndg{\boeta}}_{\mathcal{B},\lambda}+{\ndg{\boeta_t}}_{\mathcal{B},\lambda}){\ndg{\boeta^{\tilde{z}}}}_{\mathcal{B},\lambda}
\end{aligned}
\end{equation}
As for the last term, we proceed as follows.
Spelling out the definition of $\mathcal{P}$, 
we have
\begin{align*}
\langle {\mathcal{P}}|_{(t,\sbrho(t))}, {\tilde{\bz}} \rangle
&=\Bigl\langle \int_0^1{\mathcal N}'(\bw^{\theta})(\brho)\,d\theta,\tilde{\bz}\Bigr\rangle\\
&=\int_0^1\left\langle \left(\pjump(\bw^\theta)+L
 \right)'[\brho],\jumptwo{\tilde{\bz}}\right\rangle d\theta\\
&=\int_0^1\left(\int_{\Gatr}\left(
\pjump'(\bw^\theta)[\brho])+L[\brho]\right)
:\jumptwo{\tilde{\bz}}\right) d\theta.
\end{align*}
Thus, using the assumption that the Hessian $\pjump''$ is bounded, the embedding of $H^s(\Omega^j)$ into $L^{\infty}(\Omega^j)$, $j=1,2$,
and then using~\eqref{eq:dual_elliptic_regularity}, Lemma~\ref{broken_trace}, and~\eqref{ell_proj_bound_t}, we get:

\begin{equation}
\begin{aligned}
|\langle \partial_t{\mathcal{P}}|_{(t,\sbrho(t))}, {\tilde{\bz}} \rangle|
&=
\left|\int_0^1\left( \int_{\Gatr} \bigl(
\pjump''(\bw^\theta) [\partial_t \bw^\theta, \brho]\bigr) : \jumptwo{\tilde{\bz}}\right)d\theta \right|\\
&\le \int_0^1 \left( \int_{\Gatr} 
|\pjump''(\bw^\theta)| |\partial_t \bw^\theta| |\brho|
|\jumptwo{\tilde{\bz}}|\right)d\theta\\
&\le C(\sum_{j=1}^2\ltwo{\tilde{\bz}}{\infty,\Omega^j})\int_{\Gatr}(|\partial_t \brho|+|\partial_t\bw_h|)|\brho|\\
&\le C \ltwo{\brho_t}{}
(\ltwo{\brho_t}{H^{1/2}(\Gatr)}+\ltwo{\partial_t\bw_h}{H^{1/2}(\Gatr)})\ltwo{\brho}{H^{1/2}(\Gatr)^*}\\
&\le C_{\bu} \ltwo{\brho_t}{}\ltwo{\brho}{H^{1/2}(\Gatr)^*}.
\end{aligned}
\label{eq:dt_rho_z}
\end{equation}

Hence we are left in need of an estimate of $ \ltwo{\brho}{H^{1/2}(\Gatr)^*}$. This can be obtained by the following duality argument. 

Let $\hat\bz$ be the solution \eqref{eq:weak_dual_prob} with $\mbf{\alpha} = 0$ and $\mbf{\beta} = \delta(\brho)$. Here $\delta$ is the duality map \cite[IIB, p.~860]{Zeidler} between $H^{1/2}(\Gatr)$ and $H^{1/2}(\Gatr)^*$ such that $\ltwo{\delta(\brho)}{H^{1/2}(\Gatr)}=\ltwo{\brho}{H^{1/2}(\Gatr)^*}$ and $\langle \brho,\delta(\brho)\rangle_{\Gatr}=\ltwo{\brho}{H^{1/2}(\Gatr)^*}^2$. Testing in~\eqref{eq:weak_dual_prob} with $\bv=\brho$ yields
\begin{equation}
\begin{aligned}
\ltwo{\brho}{H^{1/2}(\Gatr)^*}^2&= \langle \brho,\delta(\brho)\rangle_{\Gatr}\\
&={B}(\brho,\hat\bz)+\lambda \ltwoin{\brho}{\hat\bz}
+\langle{\mathcal{P}}\brho,{\hat\bz}\rangle,\\
&={B}(\brho,\hat\bz)+\lambda \ltwoin{\brho}{\hat\bz}
+N(\bu,\hat\bz)-N(\bw_h,\hat\bz)\\
&=
B(\brho,\boeta^{\hat z})
+\lambda \ltwoin{\brho}{\boeta^{\hat z}}
+N(\bu,\boeta^{\hat z})-N(\bw_h,\boeta^{\hat z}),
\label{eq:dual_l2_tt}
\end{aligned}
\end{equation}
with $\boeta^{\hat z}=\hat\bz-\Pi\hat\bz$. 

Now a bound for $\ltwo{\brho}{H^{1/2}(\Gatr)^*}$ can be derived by following the same steps used above to get~\eqref{rho_l2_1}, yielding:
\begin{equation}\label{dual_l2_ttt}
\begin{aligned}
	\ltwo{\brho}{H^{1/2}(\Gatr)^*}^2 &\le C(1+h^2_{\max}\lambda)^{1/2} h^{s-1}_{\max}\,
	{\ndg{\boeta}}_{\mathcal{B},\lambda}\sum_{j=1}^2\ltwo{\hat\bz}{H^s(\Omega^j)}\\
	&\le C  (1+h^2_{\max}\lambda)^{1/2} h^{s-1}_{\max}\,
	{\ndg{\boeta}}_{\mathcal{B},\lambda}
	\ltwo{\delta(\brho)}{H^{1/2}(\Gatr)},
	\end{aligned}    
\end{equation}
having used~\eqref{eq:dual_elliptic_regularity} once again. Finally, by inserting in~\eqref{dual_l2_ttt}  the definition of $\delta(\brho)$ we conclude that
\begin{equation}\label{dual_l2_tttt}
	\ltwo{\brho}{H^{1/2}(\Gatr)^*}\le
	C  (1+h^2_{\max}\lambda)^{1/2} h^{s-1}_{\max}\,
	{\ndg{\boeta}}_{\mathcal{B},\lambda}.
\end{equation}

Using~\eqref{dual_l2_tttt} in~\eqref{eq:dt_rho_z} we bound last term on the right-hand side of~\eqref{eq:dual_l2_t}. Now, recalling~\eqref{eq:Ntbound}, the bound~\eqref{ell_proj_l2_bound_t} easily follows from~\eqref{eq:dual_l2_t}. 
\qed


\section{DG method for the parabolic system and its error analysis}\label{err_par_prob}
The above discussion motivates the introduction of the following IPDG-in-space method for the system (\ref{pde}), (\ref{modelbc}), (\ref{bcvec}), and (\ref{transmission}).

 For $t=0$, let $\bu_h(0)=\bw_h(0)$. For $t\in(0,T]$, find $\bu_h\equiv \bu_h(t)\in \FESpace$ such that
\begin{equation}\label{eq:semidiscrete}
\ltwoin{(\bu_h)_t}{\bv_h}+B(\bu_h,\bv_h)+N(\bu_h,\bv_h)+ \ltwoin{\mbf{F}(\bu_h)}{\bv_h}=l(\bv_h),\texte{for all}
\bv_h\in \FESpace,
\end{equation}
where
\begin{equation}\label{eq:linear}
l(\bv_h):=\!-\!\int_{\Gamma_{\ddd}}\hspace{-1.5mm}\Big((\mbf{g}_{\ddd}\otimes\mbf{n}):(A\nabla\bv_h)+(\mathcal{X}^-G_{\ddd}\mbf{B}):(\bv_h\otimes\mbf{n})-\Sigma
\mbf{g}_{\ddd}\cdot\bv_h\Big)
+\int_{\Gamma_{\dn}}\!\mbf{g}_{\dn}\cdot\bv_h,
\end{equation}
noting that $(\mbf{q}\otimes\mbf{n}):(\mbf{v}\otimes\mbf{n})= \mbf{q}\cdot \mbf{v}$, for $\mbf{q},\mbf{v}\in\mathbb{R}^n$, and having denoted
 $G_{\ddd}:=\diag(g_{\ddd}^1,\dots,g_{\ddd}^n)$.

We shall make use of the following result to treat the nonlinear reaction term.

\begin{lemma}\label{pf_lemma}
If $d=2$, let $q\in[1,\infty)$ and if $d=3$, let $q\in[1,6]$. Then, there exists a constant $C_{\rm PF}>0$, depending only on the geometry of the subdomains $\Omega^j$, $j \in \{ 1,2 \}$ and of the Dirichlet boundary, such that,
for all $v\in \SSpace^1$, we have
\begin{equation}\label{pf}
\ltwo{v}{q}^2\le C_{\rm PF}\max\{1,\alpha_{\min}^{-1}\}\ndg{v}^2.
\end{equation}
 Moreover, for $0\le\gamma< 2$, if $d=2$, and for $0\le\gamma\le 4/3$, if $d=3$, we have
\begin{equation}\label{pf_useful}
\ltwo{v}{\gamma+2}^{\gamma+2}\le C_{\rm PF}\max\{1,\alpha_{\min}^{-1} \}\ltwo{v}{}^{\gamma}\ndg{v}^2.
\end{equation}
Finally, for $d=2$ and $\gamma=2$ (which corresponds to $q=\infty$), (\ref{pf}) and (\ref{pf_useful}) hold for $v\in\FESpace$, with $C_{\rm PF}:=C |\log(\min_{\Omega}\tth)|$ for some constant $C>0$ depending only on $\Omega^i$.
\end{lemma}

{\bf Proof.}
From the assumptions on the topology of the subdomains, we can apply Theorem 3.7 from \cite{lasis-suli:03}, to deduce
\[
\ltwo{v}{q,\Omega^j}^2\le C\max\{1,\alpha_{\min}^{-1}\}\ndg{v}^2,
\]
with the boundary contribution taken as $\partial\Omega^j\backslash \Gatr$, for $1\le q<\infty$ if $d=2$ and for $1\le q \le 6 $ if $d=3$. For the case $q=\infty$ if $d=2$, we make use of the standard inverse estimate
$\ltwo{v}{\infty,\kappa}^2\le C\log|h_{\kappa}|  \ltwo{v}{H^1(\kappa)}^2$.
Setting $p=2/\gamma$, for $1/p+1/q=1$, which gives $q =2/(2-\gamma)$, H\"older's inequality implies (\ref{pf_useful}).\qed

We are now ready to prove a bound for the difference between the elliptic projection $\bw_h$ and the dG approximation $\bu_h$.

\begin{theorem}\label{sd_theorem_one}
Consider the notation of Lemma \ref{ell_proj_lemma} and let $\gamma$ as in Lemma \ref{pf_lemma}. Assume that $\bu \in L^2(0,T;\Hn^s)\cap L^\infty(0,T\times \Omega)$, $s>3/2$, $\bu_t \in L^2(0,T;L^2(\Omega))$, and that $h_{\max}$ is small enough. Then, we have
\[
         \ltwo{\mbf{\btheta}}{L^{\infty}(0,T;L^2(\Omega))}^2+C_{coer}
         \ltwo{\mbf{\btheta}}{L^2(0,T;\SSpace)}^2\le 4\delta^2\ex^{\tilde{C}T},
\]
where  $\btheta:=\bw_h-\bu_h$, $\tilde{C}:=C(\ltwo{\bu}{\infty,[0,T]\times\Omega},\lambda)$, and
\begin{equation}\label{delta}
\delta^2(t) = \int_0^T  \tilde{C}\ltwo{\brho }{}^2+ C\Big(\ltwo{\brho }{\gamma+2}^{\gamma+2}+
\lambda^{-1}\ltwo{\brho _t}{}^2\Big) \ud t .
\end{equation}
\end{theorem}

{\bf Proof.} Without loss of generality assume that $\gamma>0$; the case $\gamma=0$ follows by a simple modification of the argument presented below.

Let $\mbf{e}:=\brho+\btheta$, with $\brho:=\bu-\bw_h$ and $\btheta:=\bw_h-\bu_h$.
We note that the continuity of $\bw_h$ in the time variable is implied by the well-posedness of the elliptic projection problem, and that of $\bu_h$ by a standard local existence argument near $0$ on the resulting system of ordinary differential equations. Hence, $\btheta$ is continuous in $[0,T]$ with $\btheta(0)=\mathbf{0}$.

Orthogonality implies
\[
\ltwoin{\partial_t \mbf{e}}{\btheta}+B(\mbf{e},\btheta)
+N(\bu,\btheta)-N(\bu_h,\btheta)+\ltwoin{\mbf{F}(\bu)-\mbf{F}(\bu_h)}{\btheta}
= 0.
\]
From (\ref{ellipitc_projection}), we deduce
\begin{equation}\label{energy_identity}
\ha \frac{\ud}{\ud t}\ltwo{\btheta}{}^2+B(\btheta,\btheta)=
 \ltwoin{\mbf{F}(\bu_h)-\mbf{F}(\bu)}{\btheta}+N(\bu_h,\btheta)-N(\bw_h,\btheta)
+\ltwoin{\lambda \brho-\brho _t}{\btheta}.
\end{equation}

Making use of the inequality
\[
\int_{\Omega} a^{\gamma+1}b\le
\frac{\gamma+1}{\gamma+2}\ltwo{a}{\gamma+2}^{\gamma+2}+\frac{1}{\gamma+2}\ltwo{b}{\gamma+2}^{\gamma+2},
\]  for $a,b>0$, the nonlinear reaction term can be bounded as follows:
\begin{equation}
\begin{aligned}
|\ltwoin{\mbf{F}(\bu)-\mbf{F}(\bu_h)}{\btheta}|
&\le C \int_{\Omega} (1 + |\bu|^{\gamma}+|\bu_h|^{\gamma})|\bu-\bu_h||\btheta|\\
&\le C \int_{\Omega} (1 + |\bu|^{\gamma})|\bu-\bu_h||\btheta|+ C \int_{\Omega} |\bu-\bu_h|^{\gamma+1}|\btheta|\\
&\le C(\bu) \big(\ltwo{\brho }{}^2+\ltwo{\btheta}{}^2\big)
+ C\big(\ltwo{\brho }{\gamma+2}^{\gamma+2}+\ltwo{\btheta}{\gamma+2}^{\gamma+2}\big).
\end{aligned}
\end{equation}

Also, using the regularity of $\pjump$ and
(\ref{trace_estimate}), we have
\begin{equation}\label{nterms_bound}
\begin{aligned}
|N(\bu_h,\btheta)-N(\bw_h,\btheta)|
&\le C^p_{\mathcal{B}}\sum_{j=1}^2\ltwo{\btheta|_{\Omega^j}}{\Gatr}^2
\le \frac{1}{4}C_{\rm coer}\ndg{\btheta}^2+\frac{\lambda}{2}\ltwo{\btheta}{}^2,
\end{aligned}
\end{equation}
choosing $\epsilon$ and $\lambda$ as in the proof of Lemma \ref{ell_proj_lemma}.

We also have
\begin{equation}\label{rho_and_rhot_bound}
|\ltwoin{\lambda\brho -\brho _t}{\btheta}| \le \frac{\lambda}{2} \ltwo{\brho}{}^2+\frac{1}{2\lambda}\ltwo{\brho _t}{}^2+\lambda\ltwo{\btheta}{}^2.
\end{equation}
Lemma \ref{pf_lemma} implies
\begin{equation}\label{imbedding_sd}
\ltwo{\btheta}{\gamma+2}^{\gamma+2}
\le C_{\rm PF}\max\{1,\alpha_{\min}^{-1}\} \ltwo{\btheta}{}^{\gamma} \ndg{\btheta}^2.
\end{equation}
Using these bounds in (\ref{energy_identity}), along with the coercivity and continuity bounds from Lemmas \ref{coercivity} and \ref{continuity}, integrating the resulting inequality with respect to $t$ (and multiplying by $2$) between $0$ and $\tau$, we arrive at
\begin{equation}\label{energy_bound_indermediate_two}
\begin{aligned}
& \ltwo{\btheta(\tau)}{}^2+C_{coer}\int_0^{\tau}\ndg{\btheta}^2\\ \le&
\delta^2+\tilde{C}\int_0^{\tau}\ltwo{\btheta}{}^2
+\hat{C} \esssup_{t\in[0,\tau]}\ltwo{\btheta(t)}{}^{\gamma}\int_0^{\tau} \ndg{\btheta}^2\\
\le& \delta^2+\tilde{C}\!\int_0^{\tau}\ltwo{\btheta}{}^2
+\hat{C} \Big(\esssup_{t\in[0,\tau]}\ltwo{\btheta(t)}{}^2+C_{coer}\!\int_0^{\tau} \ndg{\btheta}^2\Big)^{1+\gamma/2},
\end{aligned}
\end{equation}
for $\delta^2$ and $\tilde{C}$ as in the statement of the theorem,
and $\hat{C}:=C\max\{1,(C_{coer}\alpha_{\min})^{-1}\}$.
From the assumption that the mesh-size $h_{\max}$ is small enough, we can have
\[
{\delta}\le (4\hat{C}\ex^{\tilde{C}T})^{-\frac{2+\gamma}{2\gamma}},
\]
which implies $(4\hat{C}\delta^2\ex^{\tilde{C}T})^{1+\gamma/2}\le \delta^2$. We now consider the set
\[
I:=\big\{\tau\in[0,T]: \esssup_{t\in[0,\tau]}\ltwo{\btheta(t)}{}^2+C_{coer}\int_0^{\tau} \ndg{\btheta}^2\le 4\delta^2 \ex^{\tilde{C}T}\big\},
\]
which is non-empty and closed due to the continuity of $\btheta$ with respect to the time variable, as $\btheta(0)=0$. We set $\tau^*=\max I$ and we suppose that
$\tau^*<T$.
Hence, for $\tau\le \tau^*$
, we have
\[
         \esssup_{t\in[0,\tau]}\ltwo{\btheta(t)}{}^2+C_{coer}\int_0^{\tau} \ndg{\btheta}^2\le 2\delta^2
   +\tilde{C}\int_0^{\tau}\ltwo{\btheta}{}^2.
\]
Gronwall's Lemma then implies
\begin{equation}\label{after_gronwall}
         \ltwo{\btheta(\tau^*)}{}^2+C_{coer}\int_0^{\tau^*} \ndg{\btheta}^2\le 2\delta^2\ex^{\tilde{C}T},
\end{equation}
setting $\tau=\tau^*$, which contradicts the hypothesis $\tau^*<T$, due to the continuity of the left-hand side of (\ref{after_gronwall}). Hence, $I=[0,T]$
and the result holds.\qed

\begin{corollary}\label{sd_cor}
Let $0\le\gamma\le 2$ if $d=2$, and $0\le\gamma\le 4/3$ if $d=3$. With the assumptions of Lemma \ref{ell_proj_lemma} and of Theorem \ref{sd_theorem_one} and $\bu|_{\kappa}\in H^1(0,T;[H^{k_{\kappa}+1}(\kappa)]^n)$, $k_{\kappa}\ge 1$, $\el$, we have
\begin{equation}\label{apriori_energy}
   \ltwo{\mbf{\bu-\bu_h}}{L^{\infty}(0,T;L^2(\Omega))}^2+C_{coer}
         \ltwo{\bu-\bu_h}{L^2(0,T;\SSpace)}^2 \le C\mathcal{E}((0,T],\mbf{h},\bu,\FESpace),
\end{equation}
for $C$ independent of $\mbf{h}$ and
\begin{equation}\label{curly_E}
\mathcal{E}((0,T],\mbf{h},\bu,\FESpace):=\su \int_0^T h_{\kappa}^{2s_{\kappa}}\big(
|\bu|_{[H^{k_{\kappa}+1}(\kappa)]^n}^2
+|\bu_t|_{[H^{k_{\kappa}+1}(\kappa)]^n}^2\big),
\end{equation} 
for $s_{\kappa}=\min\{m_{\kappa},k_{\kappa}\}$. Moreover, with the additional assumptions of Lemma \ref{ell_proj_l2_lemma}, we have
\begin{equation}\label{apriori_l2}
   \ltwo{\mbf{\bu-\bu_h}}{L^{\infty}(0,T;L^2(\Omega))}^2 \le C h_{\max}^{s-1}\mathcal{E}((0,T],\mbf{h},\bu,\FESpace),
\end{equation}
for $s\in(1,2]$ dictated by the regularity of the exact solution of the dual problem \eqref{eq:weak_dual_prob}.
\end{corollary}
\begin{proof}
We begin by using (\ref{pf}) to deduce the estimates 
$
\ltwo{\brho}{\gamma+2}^{\gamma+2}\le C \ndg{\brho}^{\gamma+2}
$, $
\ltwo{\brho}{}^{}\le C \ndg{\brho}^{}
$, and $
\ltwo{\brho_t}{}^{}\le C \ndg{\brho_t}^{}
$.
 An application of Lemma \ref{ell_proj_lemma} on the resulting
norms, along with standard approximation estimates for the error of the (orthogonal) $L^2$-projection (see, e.g., \cite{houston-etal:01}) already implies \eqref{apriori_energy}.

For \eqref{apriori_l2}, we set $ \varepsilon = \min\{2,k\gamma/(s-1)\}$. Then, we have the bound
\[
\ltwo{\brho}{\gamma+2}^{\gamma+2}\le \ltwo{\brho}{}^{2-\varepsilon}\ltwo{\brho}{2(\gamma/\varepsilon+1)}^{\gamma+\varepsilon}\le C \ltwo{\brho}{}^{2-\varepsilon}\ndg{\brho}^{\gamma+\varepsilon},
\]
which implies (for this choice of $\varepsilon$, since then $\gamma/\varepsilon+1\le 2$) that
\[
\ltwo{\brho}{\gamma+2}^{\gamma+2}\le C h^{s-1}_{\max}\mathcal{E}((0,T],\mbf{h},\bu,\FESpace).
\]
Lemma \ref{ell_proj_l2_lemma} and standard approximation estimates now imply \eqref{apriori_l2}.
\end{proof}

As it can be seen from the above proof, the use of the novel elliptic projection, introduced and analyzed in Section \ref{err_ell_problem}, plays crucial role in the treatment of the interface nonlinearity $\mbf{P}(\cdot)$, while the continuation argument, along with the Gagliardo-Nirenberg-type inequality \eqref{imbedding_sd}, treats the nonlinear reaction $\mbf{F}(\cdot)$. 

\begin{remark}
We note that the bounds in Theorem \ref{sd_theorem_one} and, correspondingly, in Corollary \ref{sd_cor} do not require any global  mesh quasi-uniformity assumptions.  This result also holds on domains without internal interfaces as in the setting of \cite{Lasis07}. There, a different continuation argument is used for deriving a priori bounds which requires the mesh to be globally quasiuniform (for non-symmetric spatial operators), albeit for a larger range of $\gamma$ than the one considered in the present work.
\end{remark}

\begin{remark}
It is interesting to note that, due to the careful use of $L^2$-projection operators in conjunction with assumption \eqref{wind_assumption}, the constant in \eqref{apriori_energy} depends  on the Pecl\'et number only through the control of the nonlinear interface conditions. Thus, the present bound produces  error control that is Pecl\'et number independent for the dG method applied to a single domain  problem. 
\end{remark}

We remark that it is possible to show optimal error estimates with less restrictive assumptions than \eqref{growth} on the growth of the reaction term. Indeed, assuming only that $F$ is locally Lipschitz, an optimal a priori bound can be proven, subject to certain conditions on the mesh.  This argument, motivated by ideas presented in~\cite{MR1701828,MR2299767}, for different problems will be considered elsewhere.

\section{Numerical examples}\label{sec:numex}

In view of the numerical tests, we introduce a fully discrete discretisation
for the system (\ref{pde}), (\ref{modelbc}), (\ref{transmission}).
The time discretization is based on a second-order linearly implicit method 
analysed in~\cite{MR1701828,Akrivis04}.

Let us write the semidiscrete DG formulation~\eqref{eq:semidiscrete} in matrix notation as
\begin{equation}\label{eq:algebric}
{\rm M \rm U_t = LU + F(U)}
\end{equation}
where ${\rm L}$ collects all linear terms and ${\rm F}$ the nonlinear terms in \eqref{eq:semidiscrete}.
Treating the linear terms with the second-order Adams-Moulton method (trapezium rule) and the
nonlinear terms with the second-order Adams-Bashforth method  yields the following
fully discrete method (AB2-AM2):
\[
{\rm M} U^{n+1}= {\rm M}U^n + k (\theta {\rm L} U^{n+1}+(1-\theta) {\rm L} U^{n})+\frac{k}{2}(3 {\rm F}(U^n)-{\rm F}(U^{n-1})).
\]

{\bf Convergence test.} We test the validity of the IPDG method and above error bounds on a system of two equations for which the exact solution is known.

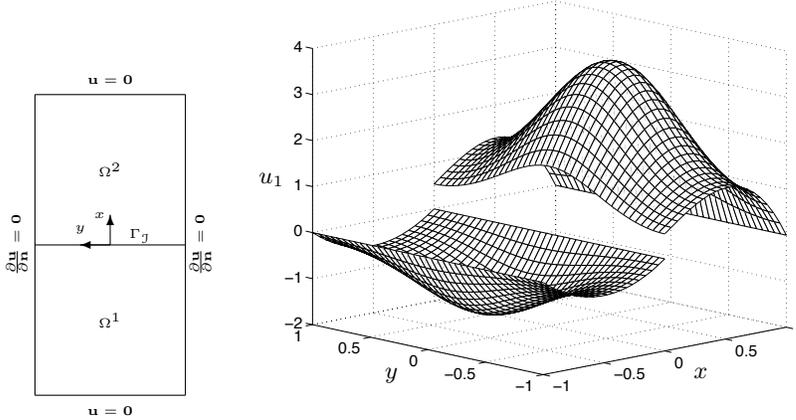
\begin{figure}[!t]
\begin{center}
\setlength{\unitlength}{2cm}
\begin{picture}(0,0)
\put(0,0){\line(0,1){2}}
\put(1,0){\line(0,1){2}}
\put(0,0){\line(1,0){1}}
\put(0,1){\line(1,0){1}}
\put(0,2){\line(1,0){1}}
\put(.5,.5){\makebox(0,0){\tiny{$\Omega^1$}}}
\put(.5,1.5){\makebox(0,0){\tiny{$\Omega^2$}}}
\put(.7,1.06){\makebox(0,0){\tiny{{{$\Gatr$}}}}}
\put(.43,1.2){\makebox(0,0){\tiny{$x$}}}
\put(.3,1.1){\makebox(0,0){\tiny{$y$}}}
\put(-.12,1){\makebox(0,0){\rotatebox{90}{\tiny{$\frac{\partial \bu}{\partial\mbf{n}}=\mbf{0}$}}}}
\put(1.1,1){\makebox(0,0){\rotatebox{90}{\tiny{$\frac{\partial \bu}{\partial\mbf{n}}=\mbf{0}$}}}}
\put(.5,1){\vector(-1,0){.2}}
\put(.5,1){\vector(0,1){.2}}
\put(.5,-.1){\makebox(0,0){\tiny{$\bu=\mbf{0}$}}}
\put(.5,2.1){\makebox(0,0){\tiny{$\bu=\mbf{0}$}}}
\end{picture}
\hspace{25 mm}
\begin{overpic}[scale=0.4]{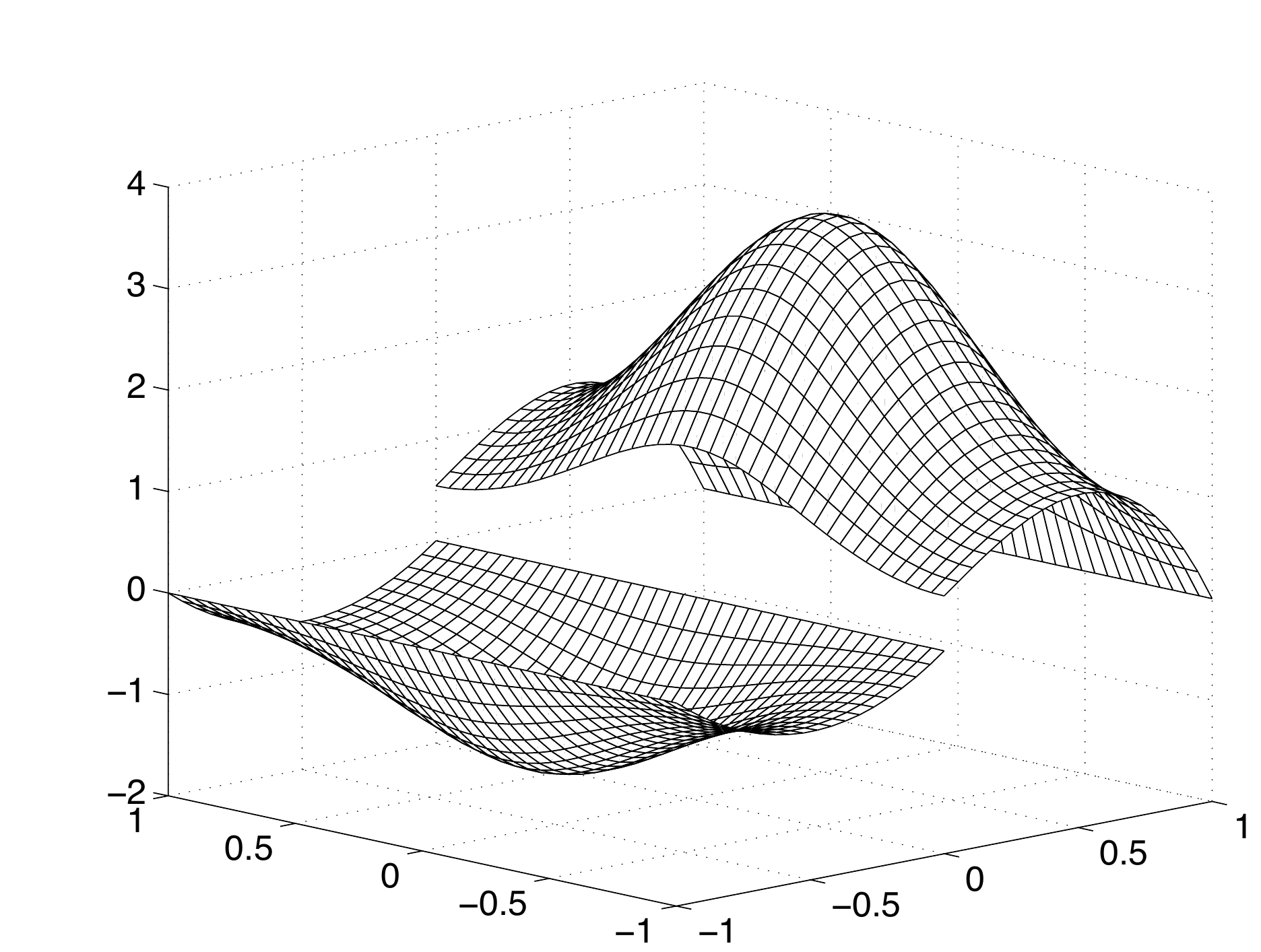}
 \small{\put(80,3){$x$}
        \put(26,3){$y$}
        \put(4,37){$u_1$}
       }
\end{overpic}
\end{center}
\caption{Convergence test. Solution domain and boundary conditions (left). First component of the exact solution at
the final time $t=1$ (right).}
\label{fig:domain}
\end{figure}

Let the domain $\Omega=[-1,1]^2$ be subdivided into two subdomains interfacing at $x=0$; thus
$\Omega^1=[-1,0]\times [-1,1]$ and $\Omega^2=[0,1]\times [-1,1]$.
We set $\Gamma_{\ddd}=\{ \pm1\}\times [-1,1]$ and 
$\Gamma_{\dn}=(-1,1)\times\{ \pm1\}$ and
impose homogeneous Dirichlet and Neumann boundary conditions
on $\Gamma_{\ddd}$ and $\Gamma_{\dn}$, respectively, as shown in Figure~\ref{fig:domain} (left).
For $t>0$ we consider the system of two advection-diffusion equations
\begin{equation}\label{ex2:system}
\left\{
\begin{array}{l}
\vspace{2 mm}
u_t - \Delta u -u_x   = f^u +
\left\{
\hspace{-2 mm}
\begin{array}{ll}
\vspace{1 mm}
  u^2 - v (1-v) &\T{in $\Omega^1$}\\
 - v &\T{in } \Omega^2,\\
\end{array}
\right.\\
v_t - \Delta v -v_x = f^v+u \quad\T{in $\Omega^1\cup \Omega^2$},
\end{array}
\right.
\end{equation}
and accordingly construct the forcing terms
$f^{u,v}:[0,1]\times\Omega\rightarrow\mathbb{R}$ in order to
yield as exact solution
\begin{equation}\label{ex2:solution}
\left(\hspace{-1 mm}
\begin{array}{c}
\vspace{1 mm}
u\\
v
\end{array}
\hspace{-1 mm}
\right)      =
\left(\hspace{-1 mm}
\begin{array}{c}
\vspace{1 mm}
\cos t\\
\sin t
\end{array}
\hspace{-1 mm}
\right)
e^{(y^2-1)^2}
\times
\left\{
\begin{array}{ll}
\vspace{1 mm}
4 x (1 + x)
 & \T{in $\Omega^1$,}\\
(- 4 x^3  + 3 x + 1) &\T{in $\Omega^2$}.
\end{array}
\right.
\end{equation}
The first component of the solution at 
time $t=1$ is shown in 
Figure~\ref{fig:domain} (right).
The functions in~\eqref{ex2:solution} are compatible
with the interface conditions~\eqref{transmission}
with respect to the interface parameters 
\[
\mbf{P}=\T{diag}(3,3), 
\qquad 
\Weight^{1} =\T{diag}(1,1),
\qquad 
\Weight^{2} =\T{diag}(0,0),
\qquad 
\Friction =\T{diag}(1,1),
\]
which are therefore used to close 
problem~\eqref{ex2:system} with 
conditions~\eqref{transmission}.
We tested the convergence rate of our method
under space discretisation refinement.
A fixed time step of size $.5\times 10^{-3}$ is used throughout while an initial uniform square $4\times 4$ mesh is uniformly refined. The value $C_{\sigma}=10$ was used throughout for the discontinuity-penalization parameter costant.
Numerical results are reported in Table~\ref{table:ex2} in the cases of bilinear and biquadratic spatial discretizations. 
The predicted error rates of convergenece are confirmed in  both the $L^2(0,1;\SSpace)$, (viz., $\|\cdot\|^2_{L^2(0,1;\SSpace)} = \int_0^1\ndg{\cdot}^2$,) and ${L^{\infty}(0,1;L^2(\Omega))}$ norms.
\begin{table}[!t]
    \label{table:ex2}
\begin{center}
\begin{tabular}{|r|r|c|c|c|c|} \hline
\# cells & \# dofs &
\multicolumn{2}{|c|}{${L^2(0,1;\SSpace)}$} &
\multicolumn{2}{|c|}{${L^{\infty}(0,1;L^2(\Omega))}$}\\ \hline
\multicolumn{6}{c}{${\tt m}\equiv 1$}\\ \hline
4 & 32 & 1.579e+01 & - & 2.498e+00 & -\\ \hline
16 & 128 & 7.617e+00 & 1.05 & 8.293e-01 & 1.59\\ \hline
64 & 512 & 3.615e+00 & 1.08 & 2.400e-01 & 1.79\\ \hline
256 & 2048 & 1.736e+00 & 1.06 & 6.497e-02 & 1.89\\ \hline
1024 & 8192 & 8.475e-01 & 1.03 & 1.693e-02 & 1.94\\ \hline
4096 & 32768 & 4.182e-01 & 1.02 & 4.324e-03 & 1.97\\ \hline
\multicolumn{6}{c}{${\tt m}\equiv 2$}\\ \hline
4 & 72 & 4.781e+00 & - & 3.812e-01 & -\\ \hline
16 & 288 & 1.013e+00 & 2.24 & 5.422e-02 & 2.81\\ \hline
64 & 1152 & 2.282e-01 & 2.15 & 7.330e-03 & 2.89\\ \hline
256 & 4608 & 5.480e-02 & 2.06 & 1.240e-03 & 2.56\\ \hline
1024 & 18432 & 1.349e-02 & 2.02 & 1.743e-04 & 2.83\\ \hline
4096 & 73728 & 3.427e-03 & 1.98 & 2.259e-05 & 2.95\\ \hline
\end{tabular}
\end{center}
\caption{Convergence test. Errors and convergence rates under uniform mesh refinement. 
A fixed time step of size $.5\times 10^{-3}$
was used.
Bilinear (above) and biquadratic (below) dG in space discretizations.
}\vspace{-1mm}
\end{table}

{\bf Advection-dominated test.} With this numerical example we test the robustness of the method in the advection-dominated regime. 
We consider once again the domain $\Omega=[-1,1]^2$ subdivided into the  subdomains $\Omega^1$ and $\Omega^2$ already considered in the previous test. In $\Omega$, we solve a scalar equation with no reaction and the 
constant diffusion $a_1=10^{-2}$ 
and the advection field $B_1=(0.5,0.5)$. 
On the boundary $\partial\Omega$, we set homogeneous Neumann boundary conditions. 
 The parameters 
\[
\mbf{P}=2/10, 
\qquad 
\Weight^{1} =5/6,
\qquad 
\Weight^{2} =1/6,
\qquad 
\Friction =6/10,
\]
are used in the  interface conditions~\eqref{transmission}. In this case, transfer across the interface is mainly advection-driven. Further, setting the friction coefficient to
less than one models the case in which the interface acts as a filtering wall on the advected quantity; hence a boundary layer in the upwind subdomain in the proximity of the interface is expected. 

\begin{figure}
\centering
\includegraphics[scale=0.12]{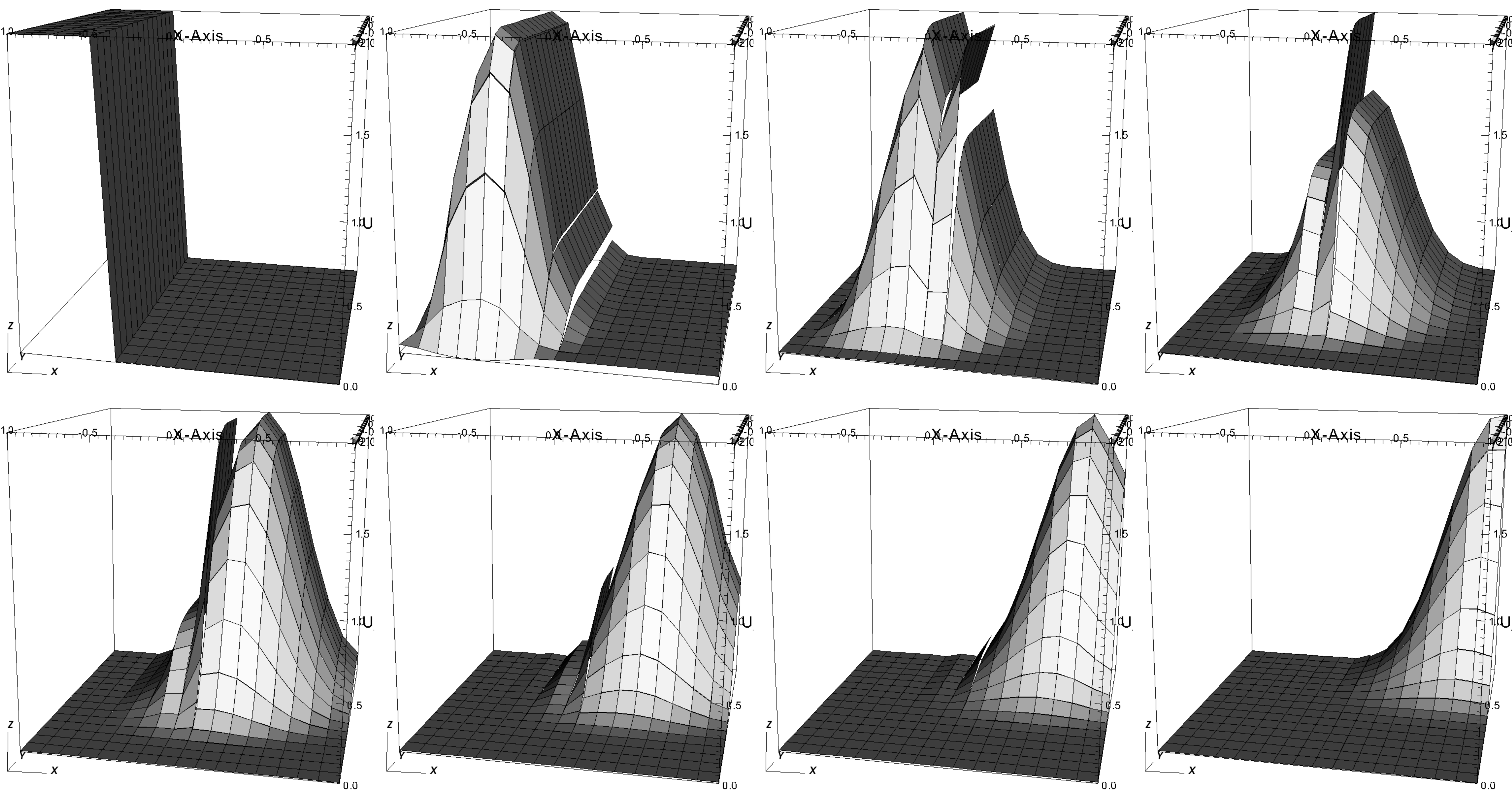}
\caption{Advection-dominated test. Snapshots of the solution computed on a uniform $16\times 16$ mesh using bilinear elements: the initial condition (top-left) followed by the solution at time intervals of $0.5$.}
\label{fig:ex2sol:16}
\end{figure}

\begin{figure}
\centering
\includegraphics[scale=0.16]{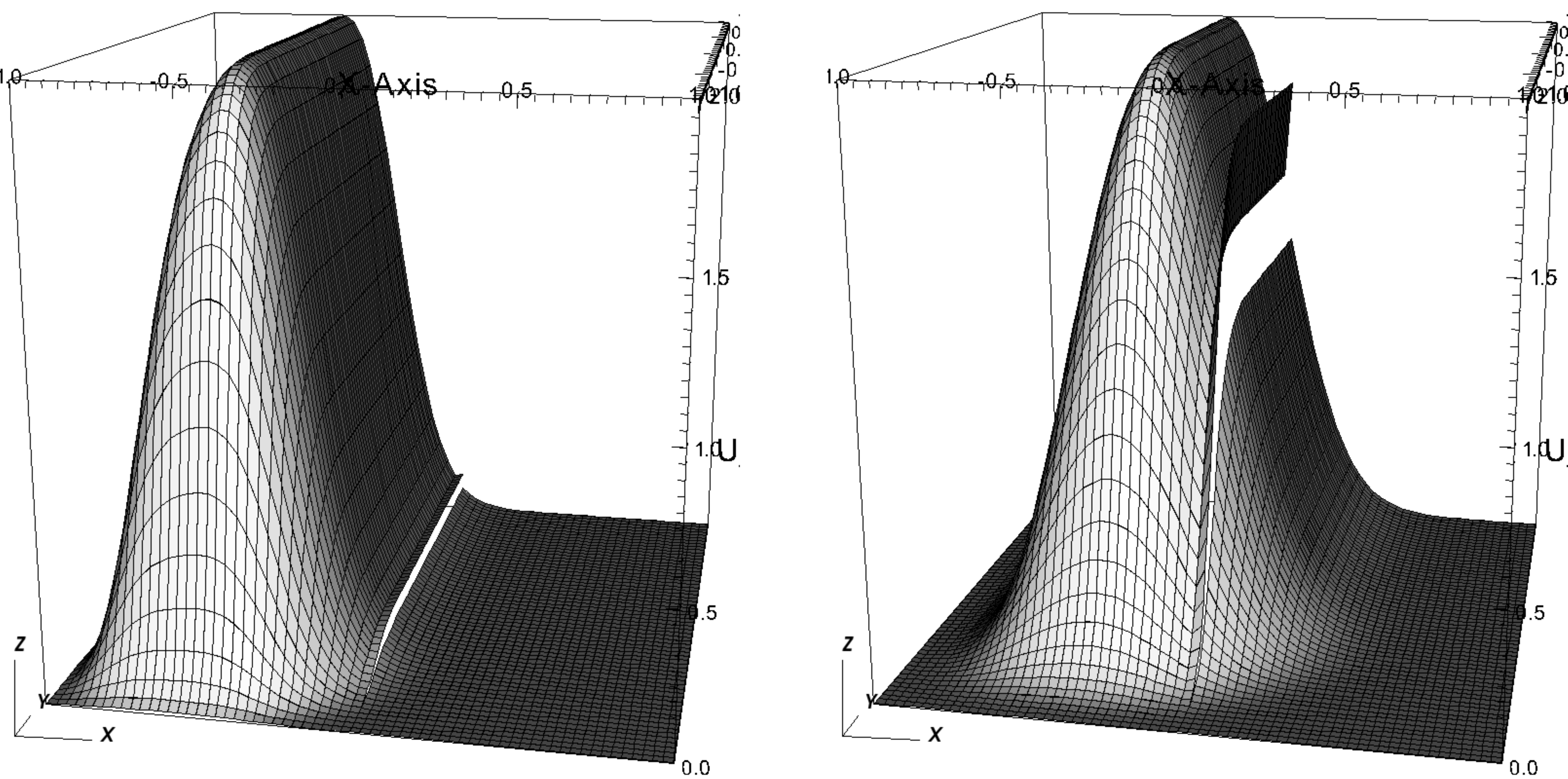}
\caption{Advection-dominated test. Snapshots of the solution computed on a uniform $64\times 64$ mesh using bilinear elements: solution at $t=0.5$ (left) and $t=1$ (right) corresponding to the second and third plots in the first row of Figure~\ref{fig:ex2sol:16}, respectively.}
\label{fig:ex2sol:64}
\end{figure}

We solve the problem on a uniform $16\times 16$ mesh using bilinear elements. Such mesh is not fine enough to resolve the layer forming in the proximity of the interface where the solution is also discontinuous, see Figure~\ref{fig:ex2sol:16}. Nevertheless, the numerical solution is stable, and the expected behaviour of the solution is accurately captured, as we can see by comparison with the solution obtained with a layer resolving $64\times 64$ mesh and shown in Figure~\ref{fig:ex2sol:64}. 
We note that the method remains stable when smaller values that $10^{-2}$  are used for the diffusion $a_1$; these results are omitted for brevity.

\section{Concluding remarks}\label{conclusions}
A dG method for the numerical solution of nonlinear interface problems modelling mass transfer through semi-permeable membranes is presented and a priori error bounds are shown under typical regularity  assumptions. The good performance of the method is highlighted  through numerical experiments. A number of extensions of the presented results can be made with modest modifications.  For instance, $hp$-version error bounds can be shown and more general convection coefficients $\mbf{B}$ can be treated. We refrained from doing so in the interest of simplicity of the presentation. 
 Interesting directions of further research are the consideration of the variational crimes due to inexact representation of the interface manifold, using, e.g., unfitted finite elements~\cite{BarEll:87,MR1622502} or the related approach in~\cite{Melenk_interface10}, and the treatment of more general interface nonlinearities. These will be considered elsewhere.

\bibliographystyle{siam}
\bibliography{bibliography}

\end{document}